\documentclass[aos,MSNbibl,rotating,seceqn,dvips]{arximspdf}

\doi{10.1214/13-AOS1135} 
\volume{41}
\issue{1}
\pubyear{2013}
\firstpage{1865}
\lastpage{1891}

\makeatletter

\newcommand{\rrVert}{\Vert}
\newcommand{\rrvert}{\vert}
\newcommand{\llVert}{\Vert}
\newcommand{\llvert}{\vert}

\newtheorem{corollary}{Corollary}
\newproclaim{defi}[prop]{Definition}
\newtheorem{theorem}{Theorem}
\newproclaim{remark}[prop]{Remark}
\newproclaim{exam}[prop]{Example}
\newtheorem{lemma}{Lemma}
\newproclaim{ass}{Assumption}
\newproclaim{example}{Example}
\makeatother

\begin{document}
\begin{frontmatter}

\title{Volatility occupation times}
\runtitle{Volatility occupation times}

\begin{aug}
\author[a]{\fnms{Jia} \snm{Li}\ead[label=e1]{jl410@duke.edu}\thanksref{t1}},
\author[b]{\fnms{Viktor} \snm{Todorov}\corref{}\thanksref{t2}\ead[label=e2]{v-todorov@northwestern.edu}}
\and
\author[a]{\fnms{George} \snm{Tauchen}\ead[label=e3]{george.tauchen@duke.edu}}
\thankstext{t1}{Supported in part by NSF Grant SES-1227448.}
\thankstext{t2}{Supported in part by NSF Grant SES-0957330.}
\runauthor{J. Li, V. Todorov and G. Tauchen}
\affiliation{Duke University, Northwestern University and Duke University}
\address[a]{J. Li\\
G. Tauchen\\
Department of Economics\\
Duke University \\
Durham, North Carolina 27708-0097\\
USA\\
\printead{e1}\\
\phantom{E-mail:\ }\printead*{e3}}

\address[b]{V. Todorov\\
Department of Finance\\
Northwestern University\\
Evanston, Illinois 60208-2001\\
USA\\
\printead{e2}}

\end{aug}

\received{\smonth{11} \syear{2012}}
\revised{\smonth{5} \syear{2013}}

%
\begin{abstract}
We propose nonparametric estimators of the occupation measure and the
occupation density of the diffusion coefficient (stochastic volatility) of a
discretely observed It\^o semimartingale on a fixed interval when the
mesh of the observation grid shrinks to zero asymptotically. In a first
step we estimate the volatility locally over blocks of shrinking
length, and then in a second step we use these estimates to construct a
sample analogue of the volatility occupation time and a kernel-based
estimator of its density. We prove the consistency of our estimators
and further derive bounds for their rates of convergence. We use these
results to estimate nonparametrically the quantiles associated with the
volatility occupation measure.
\end{abstract}

%


\end{frontmatter}

\section{Introduction}
Continuous-time It\^o semimartingales are used widely to model
stochastic processes in various areas
such as finance. The general It\^o semimartingale process is given by
%
%
\begin{equation}
\label{eqXintro} X_{t} = X_{0}+\int_{0}^{t}b_{s}\,ds+
\int_{0}^{t}\sigma_{s}\,dW_{s}+J_{t},
\end{equation}
where $b_t$ and $\sigma_t$ are
processes with c\`{a}dl\`{a}g paths, $W_t$ is a Brownian motion and
$J_t$ is a
jump process; formal
conditions are given in the next section. Inference for model~(\ref
{eqXintro}) in the
general case (either in a parametric or a nonparametric context) is
quite complicated because of the
many ``layers of latency,'' for example, as typical in financial
applications, $\sigma_t$ and $J_t$ can
have randomness not captured by~$X_t$.

When $X$ is sampled at discrete times but with the mesh of the
observation grid shrinking to zero, that is,
high-frequency data of $X$ are available, the distinct pathwise
behavior of the components in (\ref{eqXintro}) can be used to
nonparametrically separate them. Indeed,
various techniques have been already proposed to estimate
nonparametrically the integrated
variance $\int_0^T\sigma_s^2\,ds$ over a specific interval $[0,T]$ (see, e.g.,~\cite{BNS06} and
\cite{Ma1}), and more generally integrated variance measures of the form
$\int_{0}^{T}g(\sigma_{s}^2)\,ds$, where $g(\cdot)$ is a continuous
function with polynomial growth
[and there are more smoothness requirements on $g(\cdot)$ for
determining the rate of
convergence]; see Theorems 3.4.1 and 9.4.1 in~\cite{jacodprotter2012}
and the recent work of~\cite{JR}.

This paper extends the existing literature on high-frequency
nonparametric volatility estimation by
developing a nonparametric jump-robust estimator of the
occupation time of the latent volatility process $(V_{t})_{t\geq
0}\equiv(\sigma_{t}^2)_{t\geq0}$ where
the volatility occupation time is defined by
%
%
\begin{equation}
\label{eqot} F_{t}(x)=\int_{0}^{t}1_{\{V_{s}\leq x\}
}\,ds\qquad
\forall x>0,\  t\in [ 0,T ].
\end{equation}
Evidently, the right-hand side of (\ref{eqot}) is of the form
$\int_0^tg(V_s)\,ds$ with $g(v) =
1_{\{v\leq x\}}$, which unlike earlier work is a discontinuous function.

If $F_{t} ( \cdot ) $ is absolutely continuous with respect
to the Lebesgue measure, its
derivative $f_{t} ( \cdot ) $, that is, the volatility
occupation density, is well-defined.
By the Lebesgue differentiation theorem, the occupation density can be
equivalently defined as
%
%
\begin{equation}
\label{eqod} f_{t}(x)=\lim_{\varepsilon\downarrow0}\frac{1}{2\varepsilon}
\bigl(F_{t}(x+\varepsilon)-F_{t}(x-\varepsilon) \bigr).
\end{equation}
In addition to estimating $F_t(x)$, in this paper we also develop a
consistent estimator for the volatility occupation density $f_t(x)$
using the high-frequency record of $X$.

The occupation measure of the volatility process ``summarizes'' in a
convenient way the information
regarding the volatility behavior over the given time interval. Indeed,
for any bounded (or
nonnegative) Borel function $g(\cdot)$ (see, e.g., Theorem~6.4 of
\cite
{GH80}), we have
%
%
\begin{equation}
\label{eqotheorem} \int_{0}^{t}g(V_{s})\,ds=
\int_{\mathbb{R}_{+}}g(x)f_{t}(x)\,dx=\int_{\mathbb
{R}_{+}}g(x)\,dF_{t}(x).
\end{equation}
Thus, the occupation time and its density can be considered as the
pathwise analogues of the
cumulative distribution function and its density. 

Our interest in occupation times stems from the fact that they are
natural measures of risk, particularly in nonstationary settings where
invariant distributions do not exist; see, for example, the discussion
in~\cite{BP03}. Indeed, there has been a significant interest (both
theoretically and in practice) in pricing options based on the
occupation times of an underlying asset; see for example, \cite
{Dassios} and~\cite{Yor} and references therein. Here, we show how to
measure nonparametrically\vadjust{\goodbreak} occupation times associated with the
volatility risk of the price process. As a by-product, we also estimate
the corresponding quantiles of the actual path of the volatility
process over the fixed time interval. Since the pathwise volatility
quantiles are preserved under monotone transformations, they provide a
convenient way of studying the variability of the volatility and the
relationship of the latter with the volatility process
itself.\looseness=-1

We summarize our estimation procedure as follows. We first split the
fixed time interval into blocks of decreasing length and form local
estimates of the unobserved stochastic variance over each of the
blocks. The volatility estimates over the blocks are truncated
variations (see, e.g.,~\cite{Ma1} and~\cite{jacodprotter2012}), and we
further allow for adaptive choice of the truncation level that makes
use of some preliminary estimates of the stochastic variance. Then, our
estimator of the volatility occupation time is simply the empirical
cumulative distribution function of the local volatility estimates over
the blocks. Analogously, we estimate the volatility occupation density
from the local volatility estimates using kernel smoothing.

Our estimation problem can be compared with the recent work of \cite
{JR}. Jacod and Rosenbaum~\cite{JR}~show that an estimator of $\int_0^Tg(\sigma_s^2)\,ds$,
for $g(\cdot)$ a $C^3$ function, formed by plugging in local variance
estimates formed over blocks of decreasing length, can achieve the
efficient $\Delta_n^{-1/2}$ rate of convergence (for $\Delta_n$ being
the length of the high-frequency intervals). Similar to~\cite{JR},
our estimator here is formed by plugging local variance estimates in
our function of interest.

The main difference between the current work and~\cite{JR} is that in
our case the function $g(\cdot)$ in (\ref{eqotheorem}) is
discontinuous. As a result, the precision of estimating the volatility
occupation time depends on the uniform rate of recovering the
volatility process outside of the times of the ``big'' volatility jumps
(with the size of the ``big'' jumps shrinking asymptotically to zero).
Therefore, in the basic case when $X$ and $V$ are continuous, the rate
of convergence of the volatility occupation time estimator is (almost)
$\Delta_n^{-1/4}$ which, as we show in the paper, is the optimal
uniform rate for recovering the volatility trajectory from
high-frequency observations. By contrast,~\cite{JR} derive a central
limit theorem for the convergence of their estimator to $\int_0^Tg(\sigma_s^2)\,ds$ by making use of the assumed smoothness of $g$ and
applying second-order Taylor expansion of the function $g$ evaluated at
the local volatility estimator in their bias-correction and asymptotic
negligibility arguments.


Finally, our inference for the volatility occupation time and its
density can be compared with the estimation of occupation time and
density of a recurrent Markov diffusion process from discrete
observations of the process; see, for example,~\cite{FZ93} and \cite
{BP03}. The main difference is that here the state vector, and
therefore the stochastic volatility, is not fully observed. Hence, we
first need to recover nonparametrically the unobserved volatility
trajectory, and the error associated with recovering the volatility
trajectory determines the asymptotic behavior of our estimators.\vadjust{\goodbreak}

The paper is organized as follows. In Section~\ref{secsetup} we
introduce the formal setup and state our assumptions. In Section~\ref{sec-c} we develop our estimator of the volatility occupation measure
and prove its consistency. In Section~\ref{sec-r} we derive bounds for
the rate of convergence of the volatility occupation time estimator.
Section~\ref{sec-ker} derives a consistent estimator for the volatility
occupation density. Section~\ref{secmc} reports results from a Monte
Carlo study of our estimation technique. Section~\ref{secconcl}
concludes. Section~\ref{sec-pf} contains all proofs.

\section{Setup and assumptions}\label{secsetup}
We start with introducing the formal setup and stating our assumptions
about $X$. The process $X$ in (\ref{eqXintro}) is defined on a
filtered space $(\Omega,\mathcal{F},(\mathcal{F}%
_{t})_{t\geq0},\mathbb{P})$ with $b_t$ and $\sigma_t$ being adapted to
the filtration. Further, the jump component $J_t$ is defined as
%
%
\begin{eqnarray}\label{eqj}
J_{t}&=&\int
_{0}^{t}\int_{\mathbb{R}}\delta (
s,z ) 1_{ \{
\llvert \delta ( s,z ) \rrvert \leq1 \}
} ( \mu -\nu ) ( ds,dz )
\nonumber
\\[-8pt]
\\[-8pt]
\nonumber
&&{}+\int_{0}^{t}\int
_{\mathbb{R}}\delta ( s,z ) 1_{ \{ \llvert \delta ( s,z ) \rrvert
>1 \}
}\mu ( ds,dz ),
\end{eqnarray}
where $\mu$ is a Poisson measure on $\mathbb{R}_{+}\times\mathbb{R}$
with compensator $\nu$ of the form $\nu ( dt,dz )
=dt\otimes
\lambda (dz ) $ for some $\sigma$-finite measure $\lambda$
on $\mathbb{R}$ and $\delta\dvtx \Omega\times\mathbb{R}_{+}\times
\mathbb
{R}\mapsto\mathbb{R}$ is a predictable function. Regularity conditions
on $X_{t}$ are collected below.

\renewcommand{\theass}{\Alph{ass}}
\setcounter{ass}{0}
\begin{ass}\label{assA}
Let $r\in [ 0,2 ] $ be a
constant. The process $X$ is an It\^{o} semimartingale given by (\ref
{eqXintro})
and (\ref{eqj}), with $b_{t}$ locally bounded and $\sigma_{t}$ c\`
{a}dl%
\`{a}g. Moreover $\llvert \delta ( \omega,t,z ) \rrvert
\wedge1\leq\Gamma_{m} ( z ) $ for all $ ( \omega,t,z )
$ with $t\leq\tau_{m} ( \omega ) $, where $ ( \tau
_{m} ) $ is a localizing sequence of stopping times, and each
function $%
\Gamma_{m}$ on $\mathbb{R}$ satisfies $\int_{\mathbb{R}}\Gamma
_{m} (
z ) ^{r}\lambda ( dz ) <\infty$.
\end{ass}

Assumption~\ref{assA} can be viewed as a regularity-type condition. The
coefficient $r$ in Assumption~\ref{assA} controls the degree of activity of the
jump component $J$ and will play an important role in the rate of
convergence of our estimator. We note that $r$ provides an upper bound
for the (generalized) Blumenthal--Getoor index of $J_t$; see, for
example, Lemma~3.2.1 in~\cite{jacodprotter2012}. We next state our
assumption for the volatility occupation time.

\begin{ass}\label{assB}
Fix $x\geq0$. We have $F_{T} (
\cdot
) $ a.s. differentiable with derivative $f_{T} ( \cdot
) $
in a neighborhood $\mathcal{N}_{x}$ containing $x$. Moreover, $\sup_{z\in
\mathcal{N}_{x}}\mathbb{E}[f_{T} ( z ) ]<\infty$.
\end{ass}

Assumption~\ref{assB} is mainly concerned with the pathwise smoothness of the
occupation time. The differentiability condition amounts to the
existence of
occupation density $f_{T} ( \cdot ) $, which is not a strong
requirement; see, for example,~\cite{GH80} and~\cite{Protter}. The
condition $\sup_{z\in\mathcal{N}_{x}}\mathbb{E}%
[f_{T} ( z ) ]<\infty$ only requires the temporal average
(over $%
[ 0,T ] $) of the probability density of $V_{t}$ uniformly bounded
in the neighborhood~$\mathcal{N}_{x}$, which is satisfied by most
stochastic volatility
models. This condition is of course much weaker
than requiring $\mathbb{E[}\sup_{z\in\mathcal{N}_{x}}f_{T} (
z )
]<\infty$, as the latter would demand more on the pathwise regularity of
the occupation density.

\begin{example*}
Let $V_t$ solve the following stochastic
differential equation:
\[
dV_t = a_t^{V}\,dt+s(V_t)\,dW_t^{V}+dJ_t^{V},
\]
where $a_t^{V}$ is a locally bounded process, $W_t^{V}$ is a Brownian
motion, $J_t^{V}$ is a finite-variational jump process and $s
(\cdot
)$ has twice continuously differentiable reciprocal. Assume
further that $V$ has an invariant distribution which is $C^1$ in a
neighborhood of $x$. Then Assumption~\ref{assB} holds. This follows from an
application of It\^o's formula and Theorem IV.75 in~\cite{Protter}.
This example includes many parametric models of interest like the
square-root diffusion model and the more general constant elasticity of
variance model.
\end{example*}

For some of the results we will need a stronger condition on the
volatility occupation density, mainly its continuity which we state
formally in the next assumption.

\renewcommand{\theass}{\Alph{ass}$^{\prime}$}
\setcounter{ass}{1}
\begin{ass}\label{assBprime}
$F_{T} ( \cdot ) $ is a.s.
continuously differentiable on $\mathbb{R}$ with derivative
$f_{T}
(\cdot ) $.
\end{ass}

Assumption~\ref{assBprime} is harder to verify than Assumption~\ref{assB}. Necessary and
sufficient conditions for the continuity of the occupation density
(local time) of a Borel right Markov process are discussed in~\cite{EK07}.

We finally state a slightly stronger condition on the volatility
process that we will need for deriving the rate of convergence of our estimator.

\renewcommand{\theass}{\Alph{ass}}
\setcounter{ass}{2}
\begin{ass}\label{assC}
The process $\sigma_{t}$ is an It\^{o}
semimartingale with the form
\[
\sigma_{t}=\sigma_{0}+\int_{0}^{t}
\tilde{b}_{s}\,ds+\int_{0}^{t}\tilde {
\sigma%
}_{s}\,dW_{s}+\int_{0}^{t}
\tilde{\sigma}_{s}^{\prime}\,dW_{s}^{\prime
}+\int
_{0}^{t}\int_{\mathbb{R}}\tilde{
\delta} ( s,z ) ( \mu -\nu ) ( ds,dz ),
\]
where the processes $\tilde{b}$, $\tilde{\sigma}$, $\tilde{\sigma
}^{\prime}$ are locally bounded and adapted,
$W^{\prime}$ is a Brownian motion orthogonal to $W$, and $\tilde{%
\delta} ( \cdot ) $ is a predictable function. Moreover
$|\tilde{%
\delta} ( \omega,t,z ) |\wedge1\leq\tilde{\Gamma
}_{m} (
z ) $ for all $ ( \omega,t,z ) $ with $t\leq\tau
_{m} (
\omega ) $, where $ ( \tau_{m} ) $ is a localizing sequence
of stopping times, and for some $\tilde{r}\in(0,2]$, each function
$\tilde{%
\Gamma}_{m}$ on $\mathbb{R}$ satisfies $\int_{\mathbb{R}}\tilde
{\Gamma}%
_{m} ( z ) ^{\tilde{r}}\lambda ( dz ) <\infty$.
\end{ass}

Assumption~\ref{assC} assumes that $\sigma_t$ is an It\^o semimartingale, an
assumption that is satisfied by most stochastic volatility models. We
impose no restriction on the activity of the volatility jumps as well
as the dependence between $\sigma_t$ and $X_t$, a~generality that is
important in practical applications (particularly in finance).

\section{The estimator and its consistency}
\label{sec-c}
We next introduce our estimator of the volatility occupation time and
derive its consistency. We suppose that the process $X_t$ is observed
at discrete times $i\Delta_n$, $i=0,1,\ldots,$ on $[0,T]$ for a fixed
$T>0$ with the time lag $\Delta_n\rightarrow0$ when $n\rightarrow
\infty
$. The assumption for equidistant observations is merely for
simplicity, and the theoretical results that follow (except
Theorem~\ref
{thmevt}) will continue to hold in the case of irregular (but
nonrandom) sampling with $\Delta_n$ replaced by the mesh of the
irregular observation grid. In what follows the high-frequency
increment of any process $Y$ is denoted as $\Delta_{i}^{n}Y=Y_{i\Delta
_{n}}-Y_{ ( i-1 ) \Delta_{n}}$.

Our strategy of estimating $F_{T} ( \cdot ) $ is to first form
an approximation of the volatility trajectory and then use the latter
to form a sample analogue of $F_{T} ( \cdot ) $. To recover
the volatility trajectory we construct local approximations for the
spot variance process $V$ over blocks of shrinking length. To this end,
let $k_{n}$ be a sequence of integers with $k_{n}\rightarrow\infty$
and $k_{n}\Delta_{n}\rightarrow0$. Henceforth we use the shorthand
notation $u_{n}=k_{n}\Delta_{n}$. We also set a truncation process
$v_{n,t}$ verifying the following assumption, which is maintained
throughout the paper without further mention.

\begin{ass}\label{assD}
We have $v_{n,t}=\alpha_{n,t}\Delta
_{n}^{\varpi}$, where $\varpi\in ( 0,1/2 ) $ is constant,
and $%
\alpha_{n,t}$ is a strictly positive real-valued process such that for some
localizing sequence of stopping times $ ( \tau_{m} )$, $%
(\sup_{t\in [ 0,T ] }(\alpha_{n,t\wedge\tau_{m}}\vee
\alpha
_{n,t\wedge\tau_{m}}^{-1}))_{n\geq1}$ is tight for each $m\geq1$.

With this notation, for each $i=0,\ldots, \lfloor T/\Delta
_{n} \rfloor-k_{n}$, we set
%
%
\begin{equation}\label{eqvh}
\cases{ %
\displaystyle\widehat{V}_{i\Delta_{n}}^{\ast}=
\frac{1}{u_{n}}\sum_{j=1}^{k_{n}} \bigl(
\Delta_{i+j}^{n}X \bigr) ^{2},\vspace*{2pt}\cr
\displaystyle\widehat{V}_{i\Delta_{n}}=\frac{1}{%
 u_{n}}\sum_{j=1}^{k_{n}}
\bigl( \Delta_{i+j}^{n}X \bigr) ^{2}1_{ \{
\llvert \Delta_{i+j}^{n}X\rrvert \leq v_{n,i\Delta_n}
\} }.}
\end{equation}
Here, $\widehat{V}_{i\Delta_{n}}^{\ast}$ is a local approximation of $
V_{i\Delta_{n}}$ when $X$ is continuous, while $\widehat{V}_{i\Delta_{n}}$
serves the same purpose but is robust to the presence of jumps in~$X$.
As a generalization to the standard truncation-based methods (see e.g.,
Chapter~9 of~\cite{jacodprotter2012}),
we allow explicitly the truncation parameter $\alpha_{n,t}$ to be
time-varying and depend on $\{X_{i\Delta_n}\}_{i=1,\ldots,\lfloor
T/\Delta
_n\rfloor}$.
For example, one convenient and commonly used choice is to set $\alpha
_{n,t} = c \overline{\sigma}_n$,
where $c$ is a constant [typically in the range $(3,5)$], and
$\overline
{\sigma}_n$ is a preliminary estimate of
the average volatility over $[0,T]$. Assumption~\ref{assD} is verified as soon as
$\overline{\sigma}_n$ and $\overline{\sigma}_n^{-1}$ are tight.
Another possibility is to make $\alpha_{n,t}$ adaptive by setting
$\alpha_{n,t} = c \hat{\sigma}_{i}$ for $t\in[(i-1)u_n, iu_n)$,
where $\hat{\sigma}_{i}$ is a preliminary estimate for the volatility
in the local
window $[(i-1)u_n, iu_n)$. For example, one may take $\hat{\sigma}_{i}$
to be
a localized version of the Bipower variation estimator of~\cite{BNS04a}:
$\hat{\sigma}_{i} = ((\pi/2)u_n^{-1}\sum_{j=1}^{k_n} |\Delta^n_{i+j}
X||\Delta^n_{i+j+1} X|)^{1/2}$.\vadjust{\goodbreak}
The tightness requirement in Assumption~\ref{assD} can be easily fulfilled by replacing
$\hat{\sigma}_{i}$ with $(\hat{\sigma}_{i}\vee(1/C))\wedge C$ for some
pre-specified
regularization constant $C\geq1$. Finally, in the above two examples
for $\alpha_{n,t}$, we can further replace the constant $c$ with a
deterministic sequence $c_n$ increasing at a logarithmic rate as
$\Delta
_n\rightarrow0$.

We will use $\widehat{V}_{i\Delta_{n}}^{\ast}$ and $\widehat
{V}_{i\Delta_{n}}$ to approximate for the volatility trajectory within
the block. That is for $0\leq i\leq
\lfloor T/u_{n} \rfloor-1$,
%
%
\begin{equation}\label{eqvp}
\cases{ %
\widehat{V}_{t}^{\ast}=
\widehat{V}_{iu_{n}}^{\ast
}\quad\mbox{and}\quad\widehat{V}%
_{t}=\widehat{V}_{iu_{n}},&\quad  $t\in\bigl\lbrack
iu_{n}, ( i+1 ) u_{n}\bigr)$,
\vspace*{2pt}\cr
\widehat{V}_{t}^{\ast}=\widehat{V}_{ (
\lfloor
T/u_{n} \rfloor-1 ) u_{n}}^{\ast}
\quad\mbox{and}\quad\widehat {V}_{t}=\widehat{V} _{ (  \lfloor T/u_{n} \rfloor-1 ) u_{n}}, &\quad $\lfloor
T/u_{n} \rfloor u_{n}\leq t\leq T$. }\hspace*{-35pt}
\end{equation}
\end{ass}

\begin{remark}
We can alternatively define local estimators of volatility for each
$i=1,\ldots,\lfloor T/\Delta_n\rfloor$ by averaging the $k_n$ past squared
increments below the threshold. All the results in the paper, except
for Theorem~\ref{thmevt} below, will hold for this alternative way of
recovering the spot volatility.
\end{remark}

Using $\widehat{V}_{t}^{\ast}$ and $\widehat{V}_{t}$, our proposed
estimators of $F_{T} ( \cdot ) $ are defined as
\[
\widehat{F}_{n,T}^{\ast} ( x ) =\int_{0}^{T}1_{ \{
\widehat{V}%
_{s}^{\ast}\leq x \} }\,ds,\qquad
\widehat{F}_{n,T} ( x ) =\int_{0}^{T}1_{ \{ \widehat{V}_{s}\leq x \} }\,ds,\qquad
x\in\mathbb{R}.
\]

We first consider the pointwise consistency of $\widehat
{F}_{n,T}^{\ast
} ( x ) $ and $\widehat{F}_{n,T} ( x ) $. As a
matter of
fact, it is not much harder to prove a more general result as follows.

\begin{lemma}
\label{lem-1}Let $g\dvtx \mathbb{R}_{+}\mapsto [ 0,1 ] $ be a
measurable function and $D_{g}$ be the collection of discontinuity
points of
$g$. Suppose:
\begin{enumerate}[(ii)]
\item[(i)] Assumption~\ref{assA} holds for $r=2$;

\item[(ii)] for Lebesgue a.e. $t\in%
[ 0,T ] $, $\mathbb{P} ( V_{t}\in D_{g} ) =0$.

Then we have
\begin{enumerate}
\item[(a)]
%
%
\begin{equation}
\int_{0}^{T}g(\widehat{V}_{s})\,ds
\stackrel{\mathbb {P}} {\longrightarrow}%
\int_{0}^{T}g
( V_{s} ) \,ds; \label{thm-1a}
\end{equation}

\item[(b)] if, in addition, $X$ is continuous, then (\ref{thm-1a}) also holds when
replacing $\widehat{V}$ with $\widehat{V}^{\ast}$.
\end{enumerate}
\end{enumerate}
\end{lemma}

Lemma~\ref{lem-1} extends Theorem~9.4.1 in~\cite{jacodprotter2012} by
allowing for discontinuities in the test function $g ( \cdot
) $.
For fixed $x\geq0$, the pointwise consistency of $\widehat
{F}_{n,T} (
x ) $, and $\widehat{F}_{n,T}^{\ast} ( x ) $ if $X$ is
continuous, follows immediately [with $g ( \cdot )
=1_{ \{
\cdot \leq x \} }$] provided that $\mathbb{P} (
V_{t}=x )
=0$ for Lebesgue a.e. $t\in [ 0,T ] $. The uniform
consistency of $%
\widehat{F}_{n,T}^{\ast} ( \cdot ) $ and $\widehat
{F}_{n,T} (
\cdot ) $ is available if the occupation time $F_{T} (
\cdot
) $ is a.s. continuous, as shown below.

\begin{theorem}
\label{thm-1}Suppose Assumption~\ref{assA} holds for $r=2$ and $F_{T} (
\cdot ) $ is a.s. continuous. We have:\vadjust{\goodbreak}
\begin{longlist}[(a)]
\item[(a)] $\sup_{x\in\mathbb{R}}|%
\widehat{F}_{n,T} ( x ) -F_{T} ( x ) |\stackrel
{\mathbb
{P}}{%
\longrightarrow}0$;

\item[(b)] if, in addition, $X$ is continuous, then \textup{(a)} still
holds when replacing $\widehat{F}_{n,T}$ with $\widehat
{F}_{n,T}^{\ast}$.
\end{longlist}
\end{theorem}

Analogous to the classical notion of quantile for cumulative distribution
functions, the quantile of the occupation time $F_{T} ( \cdot
) $
is naturally defined as the left-continuous functional inverse of $%
F_{T} ( \cdot ) $: for $\alpha\in ( 0,T ) $, we
set $%
Q_{T} ( \alpha ) =\inf\{x\in\mathbb{R}\dvtx F_{T} (
x )
\geq
\alpha\}$. A natural estimator for $Q_{T} ( \alpha ) $ is $
\widehat{Q}_{n,T} ( \alpha ) =\inf\{x\in\mathbb
{R}\dvtx \widehat
{F}%
_{n,T} ( x ) \geq\alpha\}$, and $\widehat{Q}_{n,T}^{\ast
} ( \alpha
) $ can be defined analogously for $\widehat{F}_{n,T}^{\ast
} (
\cdot ) $. The consistency of the quantile estimators is given by
the next corollary.

\begin{corollary}
\label{cor-Q}Suppose Assumption~\ref{assA} holds for $r=2$ and $F_{T} (
\cdot ) $ is a.s. continuous. Let $\mathcal{Q}\equiv\{\alpha
\in
( 0,T ) \dvtx Q_{T} ( \cdot ) $ is continuous at
$\alpha$
a.s.$\}$. We have for each $\alpha\in\mathcal{Q}$:
\begin{longlist}[(a)]
\item[(a)] $\widehat{Q}%
_{n,T} ( \alpha ) \stackrel{\mathbb{P}}{\longrightarrow}%
Q_{T} ( \alpha ) $;

\item[(b)]if, in addition, $X$ is continuous, then $%
\widehat{Q}_{n,T}^{\ast} ( \alpha ) \stackrel{\mathbb
{P}}{%
\longrightarrow}Q_{T} ( \alpha ) $.
\end{longlist}
\end{corollary}

\section{\texorpdfstring{Rate of convergence of $\widehat{F}_{n,T}$}
{Rate of convergence of F n,T}}\label{sec-r}
We next study the rate of convergence of $\widehat{F}_{n,T} (
x )$. We first consider in Section~\ref{sec-rc} the case when $V$
is continuous and then the general case with
discontinuous $V$ is studied in Section~\ref{sec-rd}, where the uniform
rate of convergence of $\widehat{F}_{n,T}(\cdot)$ is also considered.

\subsection{The continuous volatility case and the uniform approximation of $V$}
\label{sec-rc}
In the continuous volatility case we can link the rate of convergence
of our volatility occupation time estimators with the rate of
convergence of $\widehat{V}_{t}^{\ast}$ and $\widehat{V}_{t}$ toward
$V_t$ on the space of c\`{a}dl\`{a}g functions equipped with the
uniform norm.
We denote the latter as
\[
\eta_{n}^{\ast}=\sup_{t\in [ 0,T ] }\bigl\llvert
\widehat {V}%
_{t}^{\ast}-V_{t}\bigr
\rrvert, \qquad \eta_{n}=\sup_{t\in
[
0,T%
] }\llvert
\widehat{V}_{t}-V_{t}\rrvert.
\]
The rate of convergence\vspace*{1pt} of $\widehat{F}_{n,T}$ and $\widehat{Q}_{n,T}$
is then related with $\eta_{n}$ through Lemma~\ref{lem-fu} below; an
analogous result holds for $\widehat{F}_{n,T}^{\ast}$, $\widehat
{Q}_{n,T}^{\ast}$ and $\eta_{n}^{\ast}$, but is omitted here for brevity.

\begin{lemma}\label{lem-fu}
\textup{(a)} For any $x\geq0$ and $\alpha\in ( 0,T ) $, we have
$|\widehat{F}_{n,T} ( x ) -F_{T} ( x ) |\leq
F_{T}(x+\eta_{n})-F_{T}(x-\eta_{n})$ and $|\widehat{Q}_{n,T}(\alpha
)-Q_{T}(\alpha)|\leq\eta_{n}$.

\textup{(b)} Suppose Assumption~\ref{assB} and $\eta_{n}=O_{p} ( a_{n} ) $ for
some nonrandom sequence $a_{n}\rightarrow0$. Then $\widehat
{F}_{n,T} (x ) -F_{T} ( x ) =O_{p} (
a_{n} ) $.
\end{lemma}
In view of Lemma~\ref{lem-fu}, bounding the rate of convergence of the
occupation time estimators boils down to\vadjust{\goodbreak} establishing the asymptotic
order of magnitude of $\eta_{n}$ and~$\eta_{n}^{\ast}$. Our main
result concerning the uniform approximation of the $V$ process is given
by the following theorem.

\begin{theorem}
\label{thm-u}Suppose Assumptions~\ref{assA} and~\ref{assC} with $\Delta V_s=0$ for $s\in
[0,T]$. Let $k_{n}\asymp\Delta
_{n}^{-\gamma}$ for some $\gamma\in(r\varpi+(1\vee r)(1-2\varpi
),1)$, $%
\varpi\in((1\vee r-1)/(2 ( 1\vee r ) -r),1/2)$ and $\iota
>0$ be
arbitrarily small but fixed. We have:
\begin{longlist}[(a)]
\item[(a)] $\eta_{n}=O_{p} ( a_{n} )
$, where
%
%
\begin{equation}\label{an}
\qquad a_{n}=\cases{ %
\Delta_{n}^{\gamma-1+ ( 2-r ) \varpi}
\vee\Delta _{n}^{\gamma
/2-\iota}\vee\Delta_{n}^{ ( 1-\gamma ) /2-\iota},
&\quad  $\mbox {if } r\leq1$,
\vspace*{2pt}\cr
\Delta_{n}^{\gamma/r- ( 1-\varpi ) -\iota}\vee\Delta _{n}^{\gamma/2-\iota}
\vee\Delta_{n}^{ ( 1-\gamma )
/2-\iota}, & \quad$\mbox{if } r>1$,
\vspace*{2pt}\cr
\Delta_{n}^{\gamma/2-\iota}\vee\Delta_{n}^{ ( 1-\gamma )
/2-\iota}, &\quad
$\mbox{if }  X\mbox{ is continuous;}$}
\end{equation}

\item[(b)] if $X$ is continuous, we also have $\eta_{n}^{\ast}=O_{p}
(a_{n} ) $.
\end{longlist}
\end{theorem}
When $X$ is continuous, the terms $\Delta_{n}^{\gamma/2-\iota}$ and
$\Delta_{n}^{ ( 1-\gamma ) /2-\iota}$ capture,
respectively, the
sampling variability and the discretization bias in the approximation
of the spot variance. When $X$ is discontinuous, $a_{n}$ contains an additional
term arising from the elimination of jumps, which of course depends on
the concentration of ``small'' jumps
through $r$ (recall Assumption~\ref{assA}). The conditions on $\varpi$ and
$\gamma$ imply $a_n\to0$. In particular, when $r$ is close to $2$,
$\varpi$ and $\gamma$ need to be chosen close to $1/2$ and $1$,
respectively, to ensure that $a_{n}\rightarrow0$, rendering the rate of
convergence arbitrarily slow. Nonetheless, if $X$ is discontinuous,
$\eta_{n}$ still has the same rate of convergence as in the continuous
case, that is, $\Delta_{n}^{1/4-\iota}$, provided $r\in (
0,1/2 ) $. This rate can be achieved by setting $\varpi\in
(3/ ( 8-4r ),1/2 ) $ and $\gamma=1/2$.

Of course Theorem~\ref{thm-u} provides only a bound for the rate of
convergence of $\eta_n$ and $\eta_n^*$.
The following theorem, however, establishes the exact asymptotic
distributions of $\eta_n$ and $\eta^*_n$ in a simple
model with constant volatility.
%
\begin{theorem}\label{thmevt}
Suppose:
\begin{longlist}
\item[(i)] Assumption~\ref{assA} holds with $V_t$ constant and $b_t=0$ on $[0,T]$;

\item[(ii)] $r<1/2$ and $\varpi\in(3/(8-4r),1/2)$.

Then
%
%
\begin{equation}
\label{thmevt1} \sqrt{\log\bigl(\lfloor T/u_n\rfloor\bigr)} (
\sqrt{k_n}\eta_n - \sqrt {2}Vm_{n} )
\stackrel{\mathcal{L}} {\longrightarrow} V\times \Lambda,
\end{equation}
provided $k_n \asymp\Delta_n^{-1/2}$ and where $\Lambda$ is a random
variable with c.d.f.\break $\exp(-2\exp(-x))$, and
%
%
\begin{equation}
\label{thmevt2} m_{n} = \sqrt{2\log\bigl(\lfloor T/u_n
\rfloor\bigr)} - \frac{\log(\log
(\lfloor
T/u_n\rfloor))+\log(4\pi)}{2\sqrt{2\log(\lfloor T/u_n\rfloor)}}.
\end{equation}
If we further assume \textup{(iii)} $X_t$ is continuous, then (\ref{thmevt1})
still holds with $\eta_n$ replaced by $\eta_n^*$.
\end{longlist}
\end{theorem}

\begin{remark}
Theorem~\ref{thmevt} shows that the rates given in (\ref{an}) are
almost optimal when the jumps of $X_t$ are not very
active ($r<1/2$). To be precise, we observe that (\ref{an}) suggests
$\eta_n = O_p(\Delta_n^{1/4-\iota})$ for $\iota>0$ fixed but
arbitrarily small, while the optimal rate in Theorem~\ref{thmevt}
provides a slightly sharper bound $\eta_n = O_p(\Delta_n^{1/4}\log
(\lfloor T\Delta_n^{-1/2}\rfloor)^{1/2})$.
\end{remark}


The rate of convergence of our volatility occupation time estimators
and their quantiles is a direct corollary of Lemma~\ref{lem-fu} and
Theorem~\ref{thm-u}; the proof is
omitted for brevity.
%
\begin{corollary}
\label{cor-ratec}Let $x\geq0$ and $\alpha\in(0,T)$. Suppose
Assumption~\ref{assB}
and the same setting as in Theorem~\ref{thm-u}. Then $\widehat
{F}_{n,T} ( x ) -F_{T} ( x ) $ and $\widehat
{Q}_{n,T} (
\alpha ) -Q_{T} ( \alpha ) $ are $O_{p} (
a_{n} ) $%
. If $X$ is continuous and $k_{n}\asymp\Delta_{n}^{-1/2}$, then
$\widehat{F%
}_{n,T} ( x ) -F_{T} ( x ) $, $\widehat
{F}_{n,T}^{\ast
} ( x ) -F_{T} ( x ) $, $\widehat
{Q}_{n,T}(\alpha
)-Q_{T} ( \alpha ) $ and $\widehat{Q}_{n,T}^{\ast}(\alpha
)-Q_{T} ( \alpha ) $ are $O_{p}(\Delta_{n}^{1/4-\iota})$
for $%
\iota>0$ arbitrarily small but fixed.
\end{corollary}

We should point out that in the trivial cases when $x<\inf_{t\in
[0,T]}V_t$ or $x>\sup_{t\in[0,T]}V_t$ on a given path, the error in
recovering the occupation time will become identically zero for $n$
sufficiently high (up to taking a subsequence).

\begin{remark}
More generally, we can use Theorem~\ref{thm-u} to show in the setting
of the theorem that if $\mathcal{L}\dvtx \mathcal{D}([0,T])\rightarrow
\mathbb{R}$, where $\mathcal{D}([0,T])$ is the space of c\`{a}dl\`{a}g
functions on the interval $[0,T]$ equipped with the uniform topology,
is a continuous function, we have $\mathcal{L}(\widehat{V})\stackrel
{\mathbb{P}}{\longrightarrow}\mathcal{L}(V)$. If further $|\mathcal
{L}(f)-\mathcal{L}(g)|\leq K\sup_{s\in[0,T]}|f_s-g_s|$ for any elements
$f, g\in\mathcal{D}([0,T])$ and some positive constant $K$, then the
rate of convergence of $\mathcal{L}(\widehat{V})$ to $\mathcal{L}(V)$
is bounded by the order of magnitude of $\eta_n$ (under the conditions
of Theorem~\ref{thm-u}). An example of such a function is $\mathcal
{L}(f) = \sup_{s\in[0,T]}f_s$.
\end{remark}

\subsection{The discontinuous volatility case}
\label{sec-rd}
We now turn to the general case when the volatility process contains
jumps. When $V$ is discontinuous, bounding the rate of convergence of
the volatility occupation time estimators is much less straightforward.
Lemma~\ref{lem-fu} is still valid, however, and the uniform
approximation error $\eta_{n}$ no longer vanishes asymptotically, due
to the discontinuity in $V$. This is true even if we consider the ideal
case where $\widehat{V}_{iu_{n}}=V_{iu_{n}}$, that is, perfect
pointwise approximation is available.
Indeed, the lack of uniform approximation for a discontinuous process
with its discretized version is well known in the study of convergence
of processes.

Our strategy of bounding the rate of convergence of $\widehat
{F}_{n,T} ( x )$ and $\widehat{F}_{n,T}^* ( x )$
is to
pick out the ``big'' jumps in the $V$
process and then consider the uniform rate of approximation to the
``remainder'' process. The idea is best
illustrated in the basic case where the jumps in $V$ is finitely
active. In
this case, the volatility jumps only occur within finitely many time blocks
with the form $[iu_{n}, ( i+1 ) u_{n})$ and hence their total
effect on
the estimation is $O_{p} ( u_{n} )$. On time blocks not containing
the jumps of $V$, $V$ is continuous, so Theorem~\ref{thm-u} can be
used to
provide uniform bound. The situation becomes considerably more complicated
when $V$ has infinitely active, or even infinite variational, jumps. In
this case, one needs to compute the trade-off between picking out a smaller
number of big jumps with less accurate uniform approximation to the
remainder process, and picking out a larger number of big jumps with more
accurate uniform approximation to the remainder process. The end result of
this calculation is Theorem~\ref{thm-rd} below.

\begin{theorem}
\label{thm-rd}Suppose Assumptions~\ref{assA},~\ref{assB} and~\ref{assC}. Let $k_{n}\asymp\Delta
_{n}^{-\gamma}$ for some $\gamma\in(r\varpi+(1\vee r)(1-2\varpi
),1)$, $%
\varpi\in((1\vee r-1)/(2 ( 1\vee r ) -r),1/2)$ and $\iota
>0$ be
arbitrarily small but fixed. We have:
\begin{longlist}[(a)]
\item[(a)] $\widehat{F}_{n,T} ( x )
-F_{T} ( x ) =O_{p} ( d_{n} ) $, where
$d_{n}=a_{n}\vee
\Delta_{n}^{ ( 1-\gamma ) / ( 1+\tilde{r} )
-\iota}$
and $a_{n}$ is given by (\ref{an});

\item[(b)] if $X$ is continuous, we also have $%
\widehat{F}_{n,T}^{\ast} ( x ) -F_{T} ( x )
=O_{p} (
d_{n} ) $.
\end{longlist}
\end{theorem}

Theorem~\ref{thm-rd} establishes an upper bound for the pointwise rate
of convergence of the occupation time estimators. We remind the reader
that the rate $d_n$ depends on $r$ through $a_n$; recall (\ref{an}) and
the discussion following Theorem~\ref{thm-u}. Whether the rate is
optimal or not is an open question. The rate optimality for jump-robust
estimation of the integrated variance, that is, $\int_0^T V_s\, ds =
\int_0^{\infty} x F_T(dx)$, is studied by~\cite{jacodreiss2012}.

In order to establish a uniform bound, we invoke the stronger
Assumption~\ref{assBprime} which assumes continuity of the volatility occupation density.

\begin{theorem}
\label{thm-rdu}Consider the same setting as in Theorem~\ref{thm-rd} except
with Assumption~\ref{assBprime} replacing Assumption~\ref{assBprime}. Then \textup{(a)} and \textup{(b)} in Theorem~\ref%
{thm-rd} hold uniformly in $x\in\mathbb{R}$. Moreover, for $\alpha
\in(0,T)
$ with $f_{T}(Q_{T}(\alpha))>0$ a.s., we have $|\widehat
{Q}_{n,T}(\alpha
)-Q_{T}(\alpha)|=O_{p} ( d_{n} ) $ and, if $X$ is continuous,
$|%
\widehat{Q}_{n,T}^{\ast}(\alpha)-Q_{T}(\alpha)|=O_{p} (
d_{n} ) $.
\end{theorem}

Theorem~\ref{thm-rdu} establishes a bound for the uniform rate of
convergence of the occupation time estimators; the rate of convergence
of the quantiles follows as a consequence.

%
\section{Estimation of the volatility occupation density}
\label{sec-ker}

We now turn to estimating the volatility occupation density $f_{T}(x)$.
Clearly, the uniform convergence of the occupation time (Theorem~\ref
{thm-1}%
) does not directly lead to valid estimation for the occupation density.
While the focus of the current paper is on the occupation time, we consider
the estimation of $f_{T}(\cdot)$ theoretically complementary and
empirically relevant. Our occupation density estimator is based on the local
volatility estimates $\widehat{V}_{t}$ and $\widehat{V}_{t}^{\ast}$ and
kernel smoothing. In particular, we propose the following kernel estimator
of Nadaraya--Watson type:
\[
\widehat{f}_{n,T} ( x ) \equiv\int_{0}^{T}
\frac
{1}{h_{n}}\kappa 
%
\biggl(\frac{\widehat{V}_{s}-x}{h_{n}}%
%
\biggr)\,ds,
\]
where $h_{n}\rightarrow0$ is a bandwidth sequence and the kernel
function $%
\kappa\dvtx \mathbb{R}\mapsto\mathbb{R}_{+}$ is bounded $C^{1}$ with bounded
derivative and $\int_{\mathbb{R}}\kappa ( x ) \,dx=1$. We
can define
$\widehat{f}_{n,T}^{\ast}(x)$ similarly but with $\widehat
{V}_{s}^{\ast}$
replacing $\widehat{V}_{s}$.

Below, we consider a weight function $w\dvtx \mathbb{R}\mapsto\mathbb{R}_{+}$
with $\int_{\mathbb{R}}w ( x ) \,dx<\infty$. For generic real-valued
functions $g_{1}$ and $g_{2}$ on $\mathbb{R}$, we denote
\[
\llVert g_{1}-g_{2}\rrVert _{w}\equiv\int
_{\mathbb{R}}\bigl\llvert g_{1} ( x ) -g_{2} (
x ) \bigr\rrvert w ( x ) \,dx.
\]

\begin{theorem}
\label{thm-kd}Suppose:
\begin{longlist}[(iii)]
\item[(i)] Assumptions~\ref{assA},~\ref{assBprime} and~\ref{assC};

\item[(ii)] $r\in (
0,2 ) $ and $\varpi\in((1\vee r-1)/(2(1\vee r)-r),1/2)$;

\item[(iii)] $k_{n}\asymp\Delta_{n}^{-\gamma}$ for some $\gamma\in (
0,1 ) $;

\item[(iv)] $V_{t}^{-1}$ is locally bounded;

\item[(v)] for some $\beta\in(0,1]$ and
any compact $\mathcal{K}\subset ( 0,\infty ) $, there
exists a
constant $C_{\mathcal{K}}>0$, such that for all $x,y\in\mathcal{K}$,
$%
\mathbb{E}|f_{T}(x)-f_{T}(y)|\leq C_{\mathcal{K}}|x-y|^{\beta}$;

\item[(vi)] $%
\int_{\mathbb{R}}\kappa ( z ) \llvert  z\rrvert
^{\beta
}\,dz<\infty$.
\end{longlist}
We set
\[
\bar{a}_{n}^{\ast}\equiv\Delta_{n}^{\gamma/2}
\vee\Delta _{n}^{ (
1-\gamma ) /2},\qquad \bar{a}_{n}\equiv
\bar{a}_{n}^{\ast}\vee \Delta _{n}^{{(1-r\varpi-\theta)}/{(1\vee r)}- ( 1-2\varpi ) },
\]
where $\theta=0$ when $r\leq1$ and $\theta>0$ is arbitrarily fixed
when $%
r>1$. Then for each $x\geq0$, we have:
\begin{longlist}[(a)]
\item[(a)] $\widehat{f}_{n,T} ( x )
-f_{T} ( x ) $ and $\Vert\widehat{f}_{n,T}-f_{T}\Vert_{w}$
are $%
O_{p}(h_{n}^{-2}\bar{a}_{n}\vee h_{n}^{\beta})$;

\item[(b)] if $X$ is continuous, $%
\widehat{f}_{n,T} ( x ) -f_{n,T} ( x ) $ and
$\Vert
\widehat{f}_{n,T}-f_{T}\Vert_{w}$ are $O_{p}(h_{n}^{-2}\bar
{a}_{n}^{\ast
}\vee h_{n}^{\beta})$ and moreover, the results still hold with
$\widehat{f}%
_{n,T}(\cdot)$ replaced by $\widehat{f}_{n,T}^{\ast}(\cdot)$.
\end{longlist}
\end{theorem}

\begin{remark}
Condition (v) in Theorem~\ref{thm-kd} requires the occupation density
of $%
V_{t}$ to be H\"{o}lder continuous on compacta with exponent $\beta$
under the $L_1$-norm. We
preclude the analysis for cases in which $f_{T}(\cdot)$ is differentiable,
or in a H\"{o}lder class of higher order, because occupation densities
of semimartingales in
general do not enjoy such higher-order smoothness; recall from
Assumption~\ref{assC} that $V_t$ is a semimartingale.
For example, the occupation density of a one-dimensional Brownian
motion is H\"{o}lder continuous in $L_1$ with exponent $\beta= 1/2$;
see Exercise VI.1.32 in~\cite{revuzyor}. That being said, occupation
densities of other processes, such as certain Gaussian processes (see,
e.g., Table~2 in~\cite{GH80}), may enjoy higher-order smoothness. Such
models have rarely been studied in the analysis of high-frequency
financial data and are not directly compatible with Assumption~\ref{assC}, so we
do not pursue further results here. Notice that the rate $\bar
{a}_{n}^{\ast}$ is optimized by
setting $\gamma=1/2$, resulting in $\bar{a}_{n}^{\ast}=\Delta_{n}^{1/4}$.
Furthermore, when $X$ is continuous, the estimation error of the occupation
density is $O_{p}(\Delta_{n}^{\beta/4(2+\beta)})$, which is achieved by
setting $h_{n}\asymp\Delta_{n}^{1/4(2+\beta)}$. Not surprisingly, the
smoother the occupation density (larger $\beta$), the faster the rate of
convergence.
\end{remark}

\begin{remark}
Theorem~\ref{thm-kd}(a) implies $\widehat{f}_{n,T} ( x )
-f_{T} ( x ) \stackrel{\mathbb{P}}{\longrightarrow}0$ and
$\Vert
\widehat{f}_{n,T}-f_{T}\Vert_{w}\stackrel{\mathbb
{P}}{\longrightarrow}0$,
provided that $h_{n}\rightarrow0$ and $h_{n}^{-2}\bar
{a}_{n}\rightarrow0$.
These results can be shown directly without conditions (iv)--(vi) using a
very similar proof; the details are omitted for brevity. A similar comment
applies to Theorem~\ref{thm-kd}(b).
\end{remark}

\setcounter{equation}{0}

\section{Monte Carlo}\label{secmc}
We test the performance of our nonparametric procedures on two popular
stochastic volatility models. The first is the square-root diffusion
volatility model, given by
%
%
\begin{equation}
\label{mc1} dX_t = \sqrt{V_t}\,dW_t,\qquad
  dV_t = 0.03(1.0-V_t)\,dt+0.2\sqrt{V_t}\,dB_t,
\end{equation}
$W_t$ and $B_t$ are two independent Brownian motions. Our second model
is a jump-diffusion volatility model in which the log-volatility is a
L\'{e}vy-driven Ornstein--Uhlenbeck (OU) process, that is,
%
%
\begin{equation}
\label{mc2} dX_t = e^{V_t-1}\,dW_t,\qquad
  dV_t = -0.03V_t\,dt+\,dL_t,
\end{equation}
where $L_t$ is a L\'{e}vy martingale uniquely defined by the marginal law
of $V_t$ which in turn has a selfdecomposable distribution (see Theorem~17.4 of~\cite{SATO}) with characteristic triplet (Definition~8.2 of
\cite{SATO}) of $(0,1,\nu)$ for $\nu(dx) =
\frac{2.33e^{-2.0|x|}}{|x|^{1+0.5}}1_{\{x>0\}}\,dx$ with respect to the
identity truncation function. The mean and persistence of both
volatility specifications are calibrated realistically to observed
financial data, and the two models differ in the presence of volatility
jumps as well as in the modeling of the volatility of volatility: for
model (\ref{mc1}), the transformation $\sqrt{V_t}$ is with constant
diffusion coefficient while for (\ref{mc2}) this is the case for the
transformation $\log{V_t}$.

In the Monte Carlo we fix the time span to $T=22$ days (our unit of
time is a day), equivalent to one calendar month, and we consider
$n=80$ and $n=400$, which correspond to $5$-minute and $1$-minute,
respectively, of intraday observations of $X$ in a $6.5$-hour trading
day. We set $k_n=20$ for $n=80$ and we increase it to $k_n=40$ when
$n=400$, which, respectively, correspond to $4$ and $10$ blocks per unit
of time. We finally set the truncation process at $v_{n,t} = 3\sqrt {BV_j}\Delta_n^{0.49}$
for $t\in[j-1,j)$ and where $BV_j = \frac{\pi
}{2}\sum_{i=\lfloor(j-1)/\Delta_n\rfloor+2}^{\lfloor j/\Delta
_n\rfloor
}|\Delta_{i-1}^nX||\Delta_{i}^nX|$\vspace*{2pt} is the Bipower Variation on the unit
interval $[j-1,j)$. For each realization we compute the $25$th,
$50$th and $75$th volatility qunatiles over the interval $[0,T]$. The
results from the Monte Carlo are summarized in Table~\ref{tbmc}.
Overall, the performance of our volatility quantile estimator is
satisfactory. The highest bias arises for the square-root diffusion
volatility model when volatility was started from a high value (the
$75$th quantile of its invariant distribution). Intuitively, in this
case volatility drifts toward its unconditional mean, and this results
in its larger variation over $[0,T]$, which in turn is more difficult
to accurately disentangle from the Gaussian noise in the price process,
that is, the Brownian motion $W_t$ in~$X_t$. Consistently with our
asymptotic results, the biases and the mean absolute deviations of all
volatility quantiles shrink as we increase the sampling frequency from
$n=80$ to $n=400$ in all considered scenarios.

\begin{sidewaystable}
\tablewidth=\textwidth
\caption{Monte Carlo results}
\label{tbmc}
\begin{tabular*}{\textwidth}{@{\extracolsep{\fill}}lccccccccc@{}}
\hline
 & \multicolumn{3}{c}{$\bolds{\widehat{Q}_{T,n}(0.25)}$} &
\multicolumn{3}{c}{$\bolds{\widehat{Q}_{T,n}(0.50)}$} & \multicolumn
{3}{c}{$\bolds{\widehat{Q}_{T,n}(0.75)}$}\\[-6pt]
 & \multicolumn{3}{c}{\hrulefill} &
\multicolumn{3}{c}{\hrulefill} & \multicolumn
{3}{c@{}}{\hrulefill}\\
\textbf{Start value}& \textbf{True} & \textbf{Bias} & \textbf{MAD} & \textbf{True} &
\textbf{Bias} & \textbf{MAD} & \textbf{True} & \textbf{Bias} & \textbf{MAD}\\
\hline
\multicolumn{10}{c}{Panel A: Square-root volatility model, $n=80$}
\\
$V_0 = Q^V(0.25)$ & $0.3798$ & $-0.0536$ & $0.0547$ & $0.5394$ &
$-0.0478$ & $0.0514$ & $0.7324$ & $-0.0190$ & $0.0473$ \\
$V_0 = Q^V(0.50)$ & $0.6223$ & $-0.0916$ & $0.0929$ & $0.8170$ &
$-0.0651$ & $0.0703$ & $1.0513$ & $-0.0081$ & $0.0626$ \\
$V_0 = Q^V(0.75)$ & $0.9865$ & $-0.1516$ & $0.1525$ & $1.2359$ &
$-0.0949$ & $0.1027$ & $1.5310$ & $ 0.0110$ & $0.0911$ \\[3pt]
\multicolumn{10}{c}{Panel B: Square-root volatility model, $n=400$}
\\
$V_0 = Q^V(0.25)$ & $0.3798$ & $-0.0305$ & $0.0315$ & $0.5394$ &
$-0.0304$ & $0.0327$ & $0.7324$ & $-0.0178$ & $0.0293$ \\
$V_0 = Q^V(0.50)$ & $0.6223$ & $-0.0519$ & $0.0529$ & $0.8170$ &
$-0.0412$ & $0.0453$ & $1.0513$ & $-0.0146$ & $0.0375$ \\
$V_0 = Q^V(0.75)$ & $0.9865$ & $-0.0868$ & $0.0882$ & $1.2359$ &
$-0.0596$ & $0.0654$ & $1.5310$ & $-0.0043$ & $0.0554$ \\[3pt]
\multicolumn{10}{c}{Panel C: Log-volatility model, $n=80$} \\
$V_0 = Q^V(0.25)$ & $0.1737$ & $-0.0231$ & $0.0249$ & $0.2860$ &
$-0.0269$ & $0.0302$ & $0.4519$ & $-0.0171$ & $0.0358$ \\
$V_0 = Q^V(0.50)$ & $0.3293$ & $-0.0428$ & $0.0455$ & $0.5243$ &
$-0.0460$ & $0.0524$ & $0.8069$ & $-0.0245$ & $0.0610$ \\
$V_0 = Q^V(0.75)$ & $0.6337$ & $-0.0809$ & $0.0866$ & $0.9945$ &
$-0.0807$ & $0.0968$ & $1.5162$ & $-0.0434$ & $0.1117$ \\[6pt]
\multicolumn{10}{c}{Panel D: Log-volatility model, $n=400$} \\
$V_0 = Q^V(0.25)$ & $0.1737$ & $-0.0131$ & $0.0142$ & $0.2860$ &
$-0.0158$ & $0.0180$ & $0.4519$ & $-0.0116$ & $0.0224$ \\
$V_0 = Q^V(0.50)$ & $0.3293$ & $-0.0248$ & $0.0268$ & $0.5243$ &
$-0.0276$ & $0.0318$ & $0.8069$ & $-0.0169$ & $0.0358$ \\
$V_0 = Q^V(0.75)$ & $0.6337$ & $-0.0452$ & $0.0490$ & $0.9945$ &
$-0.0480$ & $0.0575$ & $1.5162$ & $-0.0305$ & $0.0682$ \\
\hline
\end{tabular*}
\tabnotetext[]{}{\textit{Notes}: In all simulated scenarios $T=22$, and we set $k_n = 20$
for $n=80$ and $k_n = 40$ for $n=400$. In each of the cases, the
volatility is started from a fixed point being the $25$th, $50$th and
$75$th quantile of the invariant distribution of the volatility
process, denoted correspondingly as $Q^V(0.25)$, $Q^V(0.50)$ and
$Q^V(0.75)$. The columns ``True'' report the average value (across the
Monte Carlo simulations) of the true variance quantile that is
estimated; MAD stands for mean absolute deviation around true value.
The Monte Carlo replica is $1000$.}
\end{sidewaystable}

\section{Conclusion} \label{secconcl}
In this paper we propose nonparametric estimators of the volatility
occupation time and its density from discrete observations of the
process over a fixed time interval with asymptotically shrinking mesh
of the observation grid. We derive the asymptotic properties of our
volatility occupation time estimator and further invert it to estimate
the corresponding quantiles of the volatility path over the fixed time
interval. Monte Carlo shows satisfactory performance of the proposed
estimation techniques.
\section{Proofs}\label{sec-pf}
This section contains all proofs. Throughout the proof, we denote by
$K$ a
generic constant that may change from line to line. We sometimes emphasize
its dependence on some parameter $p$ by writing $K_{p}$. As is typical in
this kind of problem, by a standard localization procedure, Assumptions~\ref{assA},
\ref{assC} and~\ref{assD} can be strengthened into the following stronger versions without
loss of generality.

\renewcommand{\theass}{S\Alph{ass}}
\setcounter{ass}{0}
\begin{ass}\label{assSA}
We have Assumption~\ref{assA}. The processes $
b_{t}$ and $\sigma_{t}$ are bounded, and for some
bounded nonnegative function $\Gamma$ on $\mathbb{R}$, $|\delta
(\omega,t,z)|\leq\Gamma(z)$ and
$\int_{\mathbb{R}}\Gamma ( z ) ^{r}\lambda (
dz )
<\infty$.
\end{ass}

\renewcommand{\theass}{S\Alph{ass}}
\setcounter{ass}{2}
\begin{ass}\label{assSC}We have Assumption~\ref{assC}. The processes $
\tilde{b}_{t}$, $\tilde{\sigma}_{t}$ and $\tilde{%
\sigma}_{t}^{\prime}$ are bounded, and for some bounded
nonnegative function $\Gamma_{\sigma}$ on $\mathbb{R}$%
, $|\tilde{\delta}(\omega,t,z)|\leq\Gamma_{\sigma}(z)$ and $\int_{\mathbb{R}}\Gamma_{\sigma} ( z ) ^{\tilde{r}}\lambda
( dz ) <\infty$.
\end{ass}

\renewcommand{\theass}{S\Alph{ass}}
\begin{ass}\label{assSD}
We have Assumption~\ref{assD}. Moreover,
$\alpha
_{n,t}$ and $\alpha_{n,t}^{-1}$ are uniformly bounded for all $n,t$.
\end{ass}

\subsection{\texorpdfstring{Proofs in Section \protect\ref{sec-c}}{Proofs in Section 3}}
\mbox{}
\begin{pf*}{Proof of Lemma~\ref{lem-1}} (a) We set $\widehat
{V}%
_{t}^{+}=\widehat{V}_{iu_{n}}$ for $t\in\lbrack ( i-1 )
u_{n},iu_{n})$. Denote the left-hand side of (\ref{thm-1a}) by $S_{n}$
and $%
T_{n}= \lfloor T/u_{n} \rfloor u_{n}$. We have
\[
S_{n}=\int_{0}^{ (  \lfloor T/u_{n} \rfloor-1 )
u_{n}}g\bigl(%
\widehat{V}_{s}^{+}\bigr)\,ds+\int_{0}^{u_{n}}g(
\widehat{V}_{s})\,ds+%
\int_{T_{n}}^{T}g(
\widehat{V}_{s})\,ds.
\]
Since $g$ is bounded,%
%
%
\begin{equation}\quad
\mathbb{E}\biggl\llvert S_{n}-\int_{0}^{T}g
( V_{s} ) \,ds\biggr\rrvert \leq Ku_{n}+\int
_{0}^{ (  \lfloor T/u_{n} \rfloor
-1 )
u_{n}}\mathbb{E}\bigl\llvert g\bigl(
\widehat{V}_{s}^{+}\bigr)-g(V_{s})\bigr\rrvert
\,ds. \label{thm-1-101}
\end{equation}

Observe that for each $s\in\lbrack0, (  \lfloor
T/u_{n} \rfloor-1 ) u_{n})$, $\widehat
{V}_{s}^{+}\stackrel
{\mathbb{%
P}}{\longrightarrow}V_{s}$. To see this, we recall from Assumption~\ref{assSD}
that $\alpha_{n,t}\in[\underline{\alpha},\overline{\alpha}]$ for
some constant $\overline{\alpha}\geq\underline{\alpha}>0$. Let
$\widehat{V}_{s}^{+}(\overline{\alpha})$
and $\widehat{V}_{s}^{+}(\underline{\alpha})$ be defined as
$\widehat
{V}_{s}^{+}$
except with $\alpha_{n,t}$ replaced, respectively, by $\overline{\alpha}$
and $\underline{\alpha}$.
By Theorem~9.3.2 in~\cite{jacodprotter2012},
the right continuity\vspace*{1pt} of $V$ and $u_{n}\rightarrow0$, $\widehat
{V}_{s}^{+}(\overline{\alpha})$
and $\widehat{V}_{s}^{+}(\underline{\alpha})$ converge in probability
to $V_s$. The claim then follows
$\widehat{V}_{s}^{+}(\underline{\alpha})\leq\widehat{V}_{s}^{+}
\leq
\widehat{V}_{s}^{+}(\overline{\alpha})$.

Hence, by condition
(ii) and bounded convergence, for Lebesgue a.e. $s\in [ 0,T ]
$, $%
\mathbb{E}|g(\widehat{V}_{s}^{+})-g(V_{s})|=\mathbb{E}|(g(\widehat
{V}%
_{s}^{+})-g(V_{s}))1_{ \{ V_{s}\notin D_{g} \} }|\rightarrow0$.
Applying bounded convergence on (\ref{thm-1-101}), we readily obtain
(\ref%
{thm-1a}). Part (b) can be shown similarly.
\end{pf*}

\begin{pf*}{Proof of Theorem~\ref{thm-1}} (a) For each $x\geq
0$, $%
F_{T} ( x ) =F_{T} ( x- ) $ a.s. by the
continuity of
$%
F_{T} ( \cdot ) $. Hence,
\[
\int_{0}^{T}\mathbb{P} ( V_{s}=x )
\,ds = \mathbb{E} \bigl[ F_{T} ( x ) -F_{T} ( x- ) \bigr] =
0.
\]
Therefore, $\mathbb{P} ( V_{s}=x ) =0$ for Lebesgue a.e.
$s\in
[ 0,T ] $. By Lemma~\ref{lem-1} with $g ( \cdot
)
=1_{ \{  \cdot \leq x \} }$, $\widehat{F}_{n,T} (
x )
\stackrel{\mathbb{P}}{\longrightarrow}F_{T} ( x ) $.
Since $%
\widehat{F}_{n,T} ( \cdot ) $ and $F_{T} ( \cdot
)
$ are
increasing, and $F_{T} ( \cdot ) $ is continuous, this convergence
also holds locally uniformly. Since $V$ is c\`{a}dl\`{a}g, $\bar
{V}\equiv
\sup_{t\in [ 0,T ] }V_{t}=O_{p} ( 1 ) $. For any
$\eta>0$%
, there exists some $M>0$, such that $\mathbb{P} ( \bar
{V}>M )
<\eta$, yielding $\mathbb{P} ( T\not=F_{T} ( M )
)
<\eta
$. Hence, for any $\varepsilon>0$,
\begin{eqnarray*}
&&\limsup_{n\rightarrow\infty}\mathbb{P} \Bigl( \sup_{x\in\mathbb
{R}%
}
\bigl\llvert \widehat{F}_{n,T} ( x ) -F_{T} ( x ) \bigr
\rrvert >\varepsilon \Bigr)
\\
&&\qquad\leq\limsup_{n\rightarrow\infty}\mathbb{P} \Bigl( \sup
_{0\leq
x\leq
M}\bigl\llvert \widehat{F}_{n,T} ( x )
-F_{T} ( x ) \bigr\rrvert >\varepsilon \Bigr)
\\
&&\qquad\quad{}+\limsup_{n\rightarrow\infty}\mathbb{P} \Bigl( \sup_{x\geq
M}
\bigl\llvert \widehat{F}_{n,T} ( x ) -F_{T} ( x ) \bigr
\rrvert >\varepsilon \Bigr)
\\
&&\qquad\leq\limsup_{n\rightarrow\infty}\mathbb{P} \bigl( \bigl\llvert
\widehat{F}%
_{n,T} ( M ) -F_{T} ( M ) \bigr\rrvert
>\varepsilon /2 \bigr) +\mathbb{P} \bigl( T-F_{T} ( M ) >\varepsilon
/2 \bigr)
\\
&&\qquad<\eta.
\end{eqnarray*}
Sending $\eta\rightarrow0$, we readily derive the assertion in part (a).
Part (b) can be proved similarly.
\end{pf*}

\begin{pf*}{Proof of Corollary~\ref{cor-Q}} By Theorem~\ref
{thm-1}, $%
\widehat{F}_{n,T} ( \cdot ) \stackrel{\mathbb
{P}}{\longrightarrow}%
F_{T} ( \cdot ) $ uniformly. By a subsequence argument, we can
then assume $\widehat{F}_{n,T} ( \cdot ) \stackrel{\mathrm{a.s.}}{%
\longrightarrow}F_{T} ( \cdot ) $ uniformly without loss. The
assertion in part (a) then follows Lemma~21.2 of \cite
{vandervaart1998}. The
proof of part~(b) is similar.
\end{pf*}

\subsection{\texorpdfstring{Proofs in Section \protect\ref{sec-rc}}{Proofs in Section 4.1}}
\mbox{}
\begin{pf*}{Proof of Lemma~\ref{lem-fu}} (a) Observe%
%
%
\begin{eqnarray}\label{lem-fu-101}
%
%
F_{T} ( x-\eta_{n} )    
%
&=&\int_{0}^{T}1_{ \{ V_{s}\leq x-\eta_{n} \} }\,ds
\nonumber
\\[-8pt]
\\[-8pt]
\nonumber
%
%
&\leq&\widehat{F}_{n,T} ( x ) \leq
\int_{0}^{T}1_{ \{
V_{s}\leq
x+\eta_{n} \} }\,ds=F_{T}
( x+\eta_{n} ).%
\end{eqnarray}
Since $F_{T} ( x-\eta_{n} ) \leq F_{T} ( x )
\leq
F_{T} ( x+\eta_{n} ) $, the first assertion in part (a) readily
follows. Now consider the quantiles. By definition, $\widehat{F}_{n,T}(
\widehat{Q}_{n,T}(\alpha))\geq\alpha$. By~(\ref{lem-fu-101}), $F_{T}(
\widehat{Q}_{n,T}(\alpha)+\eta_{n})\geq\alpha$. Therefore,
$Q_{T}(\alpha
)\leq\widehat{Q}_{n,T}(\alpha)+\eta_{n}$. For any $\varepsilon>0$,
by (%
\ref{lem-fu-101}), we have $F_{T}(\widehat{Q}_{n,T}(\alpha)-\eta
_{n}-\varepsilon)\leq\widehat{F}_{n,T}(\widehat{Q}_{n,T}(\alpha
)-\varepsilon)<\alpha$. Hence, $\widehat{Q}_{n,T}(\alpha)-\eta
_{n}-\varepsilon<Q_{T}(\alpha)$. Since $\varepsilon>0$ is arbitrary,
$%
\widehat{Q}_{n,T}(\alpha)-\eta_{n}\leq Q_{T}(\alpha)$. The second
assertion of part (a) is then obvious.

(b) Fix some $\varepsilon>0$. There exists $M>0$ such that $\mathbb{P}
( \eta_{n}\geq Ma_{n} ) <\varepsilon/2$ for $n$ sufficiently
large. Since $a_{n}\rightarrow0$, $ [ x-\eta_{n},x+\eta
_{n}
] $
is contained in $\mathcal{N}_{x}$ with probability approaching one (w.p.a.1)
and by part (a),
\[
\bigl\llvert \widehat{F}_{n,T} ( x ) -F_{T} ( x ) \bigr
\rrvert \leq\int_{x-\eta_{n}}^{x+\eta_{n}}f_{T} ( z )
\,dz.
\]
Let $M^{\prime}=4M\sup_{z\in\mathcal{N}_{x}}\mathbb{E} [
f_{T} (
z )  ] /\varepsilon$. We have for $n$ sufficiently large,
\begin{eqnarray*}
%
%
&&\mathbb{P} \biggl( \int_{x-\eta_{n}}^{x+\eta_{n}}f_{T} (
z ) \,dz>M^{\prime}a_{n} \biggr)
\\
%
&&\qquad\leq\mathbb{P} \biggl(
\int_{x-Ma_{n}}^{x+Ma_{n}}f_{T} ( z )
\,dz>M^{\prime}a_{n} \biggr) +\mathbb{P} ( \eta_{n}\geq
Ma_{n} )
\\
%
&&\qquad<\frac{2M\sup_{z\in\mathcal{N}_{x}}\mathbb{E} [
f_{T} ( z )  ] }{M^{\prime}}+
\varepsilon/2
\\
%
&&\qquad\leq\varepsilon.%
\end{eqnarray*}
Hence, $\int_{x-\eta_{n}}^{x+\eta_{n}}f_{T} ( z )
\,dz=O_{p} (
a_{n} ) $. The assertion in part (b) then readily follows.
\end{pf*}

We now prove Theorem~\ref{thm-u}, starting with two lemmas. Below,
$\Vert \cdot \Vert_p$ denotes the $L_p$ norm.

\begin{lemma}
\label{lem-vstar}Let $p\geq1$ be a constant and $k_{n}\asymp\Delta
_{n}^{-\gamma}$ for some $\gamma\in ( 0,1 ) $. Suppose
Assumption~\ref{assSA} holds with $X$ continuous and Assumption~\ref{assSD}. Then for
each $0\leq i\leq
\lfloor T/u_{n} \rfloor-1$,%
%
%
\begin{equation}
\bigl\llvert \widehat{V}_{iu_{n}}^{\ast}-V_{iu_{n}}\bigr
\rrvert \vee\llvert \widehat{V}_{iu_{n}}-V_{iu_{n}}
\rrvert \leq \xi_{n,i}+ \sup_{s\in\lbrack iu_{n}, ( i+1 )
u_{n})}\llvert
V_{s}-V_{iu_{n}}\rrvert, \label{lem-vstar-a}
\end{equation}
where the variable $\xi_{n,i}$ satisfies $\Vert\xi_{n,i} \Vert_p
\leq
K_p k_n^{-1/2}$.
If we further have Assumption~\ref{assSC} with $\Delta V_s =0$ for $s\in[0,T]$,
then the majorant side of the above can be
bounded by $K_{p}(k_{n}^{-1/2}+u_{n}^{1/2})$ in $L_p$.
\end{lemma}

\begin{pf}
By It\^o's formula, $\widehat
{V}_{iu_{n}}^{\ast}-V_{iu_{n}} =\zeta_{n,i}^{\prime}+\zeta
_{n,i}^{\prime\prime}$, where
\begin{eqnarray*}
\zeta_{n,i}^{\prime} &=&\frac{2}{u_{n}}\int
_{iu_{n}}^{ (
i+1 )
u_{n}} ( X_{s}-X_{n,s} )
\,dX_{s},
\\
\zeta_{n,i}^{\prime\prime} &=&\frac{1}{u_{n}}\int
_{iu_{n}}^{ (
i+1 ) u_{n}} ( V_{s}-V_{iu_{n}} )
\,ds,
\end{eqnarray*}
and $X_{n,s}$ is the discretized process given by
$X_{n,s}=X_{iu_{n}+ ( j-1 ) \Delta_{n}}$ when $s\in
\lbrack iu_{n}+ ( j-1 ) \Delta_{n},iu_{n}+j\Delta_{n})$.
By classical estimates (note that $X$ is continuous),
\begin{eqnarray*}
\mathbb{E}\biggl\llvert \frac{2}{u_{n}}\int_{iu_{n}}^{ (
i+1 )
u_{n}}
( X_{s}-X_{n,s} ) b_{s}\,ds\biggr\rrvert
^{p} &\leq &K_{p}\Delta _{n}^{p/2},
\\
\mathbb{E}\biggl\llvert \frac{2}{u_{n}}\int_{iu_{n}}^{ (
i+1 )
u_{n}}
( X_{s}-X_{n,s} ) \sigma_{s}\,dW_{s}
\biggr\rrvert ^{p} &\leq &K_{p}k_{n}^{-p/2}.
\end{eqnarray*}
Since $k_{n}=o ( \Delta_{n}^{-1} ) $, we have $\Vert\zeta
_{n,i}^{\prime}\Vert_{p}\leq K_{p}k_{n}^{-1/2}$ by Minkowski's
inequality. We also observe $\vert\zeta_{n,i}^{\prime\prime}\vert
\leq\sup_{s\in\lbrack iu_{n},(
i+1) u_{n})}\vert V_{s}-V_{iu_{n}}\vert$.

Now note that $\widehat{V}_{iu_{n}}^{\ast}-\widehat
{V}_{iu_{n}}=k_{n}^{-1}%
\sum_{j=1}^{k_{n}}(\Delta_{ik_{n}+j}^{n}X/\Delta
_{n}^{1/2})^{2}1_{\{|\Delta_{ik_{n}+j}^{n}X|>v_{n,iu_n}\}}$. Under
Assumption~\ref{assSD},
$\alpha_{n,t}\geq\underline{\alpha}$ for some constant $\underline
{\alpha}>0$. Hence,
$v_{n,t}\geq\underline{v}_n\equiv\underline{\alpha} \Delta
_{n}^{\varpi}$.
Since $X$ is
continuous, for any $q\geq0$,
\begin{eqnarray*}
&&\bigl\llVert \bigl(\Delta_{ik_{n}+j}^{n}X/\Delta_{n}^{1/2}
\bigr)^{2}1_{ \{
\llvert \Delta_{ik_{n}+j}^{n}X\rrvert >v_{n,iu_n} \}
}\bigr\rrVert _{p}\\
&&\qquad\leq \biggl(
\frac{\mathbb{E}\llvert \Delta
_{ik_{n}+j}^{n}X/\Delta
_{n}^{1/2}\rrvert ^{2p+q}}{ ( \underline{v}_{n}/\Delta
_{n}^{1/2} ) ^{q}}%
\biggr) ^{1/p}
\\
&&\qquad\leq K_{p,q}\Delta_{n}^{q ( 1/2-\varpi ) /p}.
\end{eqnarray*}
Since $\varpi\in ( 0,1/2 ) $, when $q$ is taken sufficiently
large, terms in the above display can be further bounded by $%
K_{p}k_{n}^{-1/2}$. Hence, $\Vert\widehat{V}_{iu_{n}}^{\ast
}-\widehat
{V}%
_{iu_{n}}\Vert_{p}\leq K_{p}k_{n}^{-1/2}$. The first assertion then readily
follows by setting $\xi_{n,i} = |\zeta^{\prime}_{n,i}| + \vert
\widehat
{V}_{iu_{n}}^{\ast}-\widehat{V}%
_{iu_{n}} \vert$.

Now, suppose Assumption~\ref{assSC} together with $V$ being continuous. By
standard estimates, for
each $p\geq1$, the second term on the right-hand side of (\ref
{lem-vstar-a}%
) can be bounded by $K_{p}u_{n}^{1/2}$ in $L_p$. The second assertion
of the lemma is
then obvious.
\end{pf}

Under Assumption~\ref{assSA}, we set
\begin{eqnarray*}
X_{t}^{\prime} &=&X_{t}-X_{t}^{\prime\prime},
\\
X_{t}^{\prime\prime} &=&\cases{ %
\displaystyle \int
_{0}^{t}\int_{\mathbb{R}}\delta (
s,z ) \mu ( ds,dz ), & \quad$\mbox{if }  r\leq1$,
\vspace*{2pt}\cr
\displaystyle\int_{0}^{t}\int_{\mathbb{R}}
\delta ( s,z ) ( \mu -\nu ) ( ds,dz ), &\quad  $\mbox{if } r>1.$}
\end{eqnarray*}
We define $\widehat{V}^{\prime}$ as $\widehat{V}^{\ast}$ in (\ref{eqvh})
but with $X^{\prime}$ in place of $X$; in particular, $\widehat{V}%
_{iu_{n}}^{\prime}\equiv u_{n}^{-1}\sum_{j=1}^{k_{n}}(\Delta
_{ik_{n}+j}^{n}X^{\prime})^{2}$.

\begin{lemma}
\label{lem-eoj}Suppose that Assumption~\ref{assSA} holds with some $r\in(0,2)$
and Assumption~\ref{assSD}. Let $%
p\geq r\vee1$ and $\varpi\in ( \frac{p-1}{2p-r},\frac
{1}{2}
) $. Let
$\theta\in ( 0,\infty ) $ be arbitrarily fixed if $r>1$ and
$%
\theta=0$ if $r\leq1$. We have for each $i$,
\[
\bigl\llVert \widehat{V}_{iu_{n}}-\widehat{V}_{iu_{n}}^{\prime}
\bigr\rrVert _{p}\leq K_{p}\Delta_{n}^{{(1-r\varpi-\theta)}/{p}- (
1-2\varpi
) }.
\]
\end{lemma}

\begin{pf}
Under Assumption~\ref{assSD}, $\alpha_{n,t}\in
\lbrack
\underline{\alpha},\overline{\alpha}]$ for constants $\overline
{\alpha}%
\geq\underline{\alpha}>0$. We set $\bar{v}_{n}=\overline{\alpha
}\Delta
_{n}^{\varpi}$ and $\underline{v}_{n}=\underline{\alpha}\Delta
_{n}^{\varpi}$. By applying Lemma~13.2.6 in~\cite{jacodprotter2012}
[with $%
s=1$, $s^{\prime}=2$, $m=p$, $p^{\prime}=1$, $k=1$, $F ( x )
=x^{2}$], we have
\begin{eqnarray*}
&&\mathbb{E} \bigl\vert \bigl( \Delta
_{ik_{n}+j}^{n}X/\Delta _{n}^{1/2} \bigr)
^{2}1_{ \{ \llvert \Delta
_{ik_{n}+j}^{n}X\rrvert \leq\bar{v}_{n} \} }
- \bigl( \Delta _{ik_{n}+j}^{n}X^{\prime}/
\Delta_{n}^{1/2} \bigr) ^{2}1_{ \{
\llvert \Delta_{ik_{n}+j}^{n}X^{\prime}\rrvert \leq\bar{v}%
_{n} \} }
\bigr\vert^{p}
\\
&&\qquad\leq K_{p}\Delta_{n}^{{(2-r)}/{2}-\theta
}+K_{p}
\Delta_{n}^{1-r\varpi
-p ( 1-2\varpi ) -\theta}
\\
&&\qquad\leq K_{p}\Delta_{n}^{1-r\varpi-p (
1-2\varpi
) -\theta}.
\end{eqnarray*}
By a similar argument as in the proof of Lemma~\ref{lem-vstar},
\[
\bigl\llVert \bigl(\Delta_{ik_{n}+j}^{n}X^{\prime}/
\Delta_{n}^{1/2}\bigr)^{2}1_{\{|\Delta
_{ik_{n}+j}^{n}X^{\prime}|>\bar{v}_{n}\}}\bigr
\rrVert _{p}\leq K_{p,q}\Delta _{n}^{q ( 1/2-\varpi ) /p}
\]
for any $q\geq0$. Taking $q$
sufficiently large, we then derive%
%
%
\begin{eqnarray}%
\label{lem-eoj-101}
%
%
&&\bigl\llVert \bigl( \Delta_{ik_{n}+j}^{n}X/\Delta_{n}^{1/2}
\bigr) ^{2}1_{ \{ \llvert \Delta_{ik_{n}+j}^{n}X\rrvert \leq
\bar
{v}%
_{n} \} }- \bigl( \Delta_{ik_{n}+j}^{n}X^{\prime}/
\Delta _{n}^{1/2} \bigr) ^{2}\bigr\rrVert
_{p}
\nonumber
\\[-8pt]
\\[-8pt]
\nonumber
%
&&\qquad\leq K_{p}
\Delta_{n}^{{(1-r\varpi)}/{p}- ( 1-2\varpi
) -%
{\theta}/{p}}.
\end{eqnarray}
By a similar argument, we can derive (\ref{lem-eoj-101}) when $\bar{v}_{n}$
is replaced with $\underline{v}_{n}$. Since $\underline{v}_{n}\leq
v_{n,iu_n}\leq
\bar{v}_{n}$, (\ref{lem-eoj-101}) also holds when $\bar{v}_{n}$ is replaced
with $v_{n,iu_n}$. The assertion of the lemma then follows from Minkowski's
inequality.\vadjust{\goodbreak}
\end{pf}

\begin{pf*}{Proof of Theorem~\ref{thm-u}} \textit{Step} 1. We first
suppose $X$
is continuous. Observe that
\[
\eta_{n}\leq\sup_{0\leq i\leq \lfloor T/u_{n} \rfloor
-1}\llvert
\widehat{V}_{iu_{n}}-V_{iu_{n}}\rrvert +2\sup_{0\leq
i\leq
\lfloor T/u_{n} \rfloor}
\sup_{s\in\lbrack iu_{n}, (
i+1 ) u_{n})}\llvert V_{s}-V_{iu_{n}}\rrvert.
\]
Since $V$ is continuous, $\Vert\sup_{s\in\lbrack iu_{n}, (
i+1 )
u_{n})}\llvert  V_{s}-V_{iu_{n}}\rrvert \Vert_{p}\leq
K_{p}u_{n}^{1/2} $ for any $p\geq1$ by standard estimates. By Lemma~\ref%
{lem-vstar} and the maximal inequality (e.g., Lemma~2.2.2 in \cite%
{vandervaartwellner}), for any $p\geq1$, $\llVert \eta_{n}\rrVert
_{p}\leq K_{p}u_{n}^{-1/p}(k_{n}^{-1/2}+u_{n}^{1/2})$. Since
$k_{n}\asymp
\Delta_{n}^{-\gamma}$ by assumption, we derive $\Vert\eta_{n}\Vert
_{p}\leq K\Delta_{n}^{ ( \gamma\wedge ( 1-\gamma )
)
/2-\iota} $ by taking $p$ sufficiently large. The same argument yields
$%
\Vert\eta_{n}^{\ast}\Vert_{p}\leq K\Delta_{n}^{ ( \gamma
\wedge
( 1-\gamma )  ) /2-\iota}$. This finishes the proof
of part
(a) with $X$ continuous, as well as part (b).

\textit{Step} 2. We now consider part (a) allowing $X$ to be discontinuous.
Let $\widehat{V}^{\prime}$ be defined as
in Lemma~\ref{lem-eoj}. Observe that
\begin{eqnarray*}
\eta_{n} &\leq&\sup_{0\leq i\leq \lfloor T/u_{n} \rfloor
-1}\bigl\llvert
\widehat{V}_{iu_{n}}-\widehat{V}_{iu_{n}}^{\prime
}\bigr\rrvert
+\eta_{n}^{\prime},\qquad \mbox{where}
\\
\eta_{n}^{\prime} &=&\sup_{0\leq i\leq \lfloor T/u_{n}
\rfloor
-1}\bigl\llvert
\widehat{V}_{iu_{n}}^{\prime}-V_{iu_{n}}\bigr\rrvert +2\sup
_{0\leq i\leq \lfloor T/u_{n} \rfloor}\sup_{s\in
\lbrack
iu_{n}, ( i+1 ) u_{n})}\llvert V_{s}-V_{iu_{n}}
\rrvert.
\end{eqnarray*}
A similar argument as in part (a) yields $\eta_{n}^{\prime}=
O_p(\Delta_{n}^{ ( \gamma\wedge ( 1-\gamma )  )
/2-\iota})$. By the maximal inequality and Lemma~\ref{lem-eoj} for
$p=1\vee
r $,
\begin{eqnarray*}
\Bigl\llVert \sup_{0\leq i\leq \lfloor T/u_{n} \rfloor
-1}\bigl\llvert \widehat{V}_{iu_{n}}-
\widehat{V}_{iu_{n}}^{\prime}\bigr\rrvert \Bigr\rrVert
_{p} &\leq&Ku_{n}^{-1/p}\Delta_{n}^{{(1-r\varpi-\theta)
}/{p}- (
1-2\varpi ) }
\\
&\leq&K\Delta_{n}^{{(\gamma-r\varpi-\theta)}/{p}- (
1-2\varpi
) }
\\
&\leq&\cases{ %
K\Delta_{n}^{\gamma-r\varpi- ( 1-2\varpi ) },
& \quad$\mbox {if } r\leq1$,
\vspace*{2pt}\cr
K\Delta_{n}^{ ( \gamma-\theta ) /r- ( 1-\varpi
)
}, & \quad$\mbox{if } r>1.$}
\end{eqnarray*}
Taking $\theta$ sufficiently small in the $r>1$ case, we readily
derive the
assertion in part~(a).
\end{pf*}

\begin{pf*}{Proof of Theorem~\ref{thmevt}}
\textit{Step} 1. We first prove the assertion on $\eta_n^*$, so condition (iii)
is in force. In the constant volatility setting of the theorem, we have
$\sqrt{k_n}\eta_n^* = \sqrt{2}V\times M_n$, where we denote
%
%
\begin{eqnarray}
M_n &=& \sup_{i=0,\ldots,\lfloor{T}/{u_n}\rfloor-1}\bigl\llvert Z_i^n
\bigr\rrvert,
\nonumber
\\[-8pt]
\\[-8pt]
\nonumber
Z_i^n &=& \frac{\sqrt{k_n}}{\sqrt{2}V} \bigl(
\widehat{V}^*_{iu_n} -V \bigr),\qquad
i=0,\ldots,\lfloor T/u_n
\rfloor-1.
\end{eqnarray}

Under our constant volatility assumption $\{Z_i^n\}_i$ are independent
and identically distributed with distribution which is approximately
standard normal. Therefore, we can use Edgeworth expansion of the c.d.f.
together with extreme value theory to pin down the limit distribution
of $M_n$. To this end, we set
%
%
\begin{equation}
\cases{ %
\displaystyle c_{n} = \bigl(2
\log(b_n)\bigr)^{-1/2},\qquad  m_{n} = \sqrt{2
\log(b_n)} - \frac
{\log
(\log(b_n))+\log(4\pi)}{2\sqrt{2\log(b_n)}},
\vspace*{2pt}\cr
\tau_{n}(x) = c_{n}x+m_{n},  \qquad b_n =
\lfloor T/u_n\rfloor,\qquad  x\in\mathbb{R}_+. }\hspace*{-35pt}
\end{equation}
Note that $\tau_{n}(x)\asymp\sqrt{2\log(b_n)}$ and hence increases to
infinity as the number of blocks increases to infinity for every fixed $x$.

Using second-order Edgeworth expansion and denoting with $\Phi(\cdot)$,
the c.d.f. of standard normal random variable (see Theorem~2.2 and Lemma~5.4 of~\cite{Hall}), we have
%
%
\begin{equation}
\mathbb{P} \bigl(\bigl|Z_i^n\bigr|\leq\tau_n(x)
\bigr) = \Phi\bigl(\tau _n(x)\bigr)-\Phi \bigl(-\tau_n(x)
\bigr) + \frac{(\log(b_n))^4}{b_n\sqrt{k_n}}K(x)+o \biggl(\frac
{1}{k_n} \biggr)\hspace*{-35pt}
\end{equation}
for any $x$ where $K(x)$ is a polynomial of $x$. Then we have
%
%
\begin{equation}
\label{eqproofevtt} \lim_{n\rightarrow\infty} \biggl[\frac{\mathbb{P}
(|Z_i^n|\leq\tau
_n(x)  )}{\Phi(\tau_n(x))-\Phi(-\tau_n(x))}
\biggr]^{b_n} = 1,
\end{equation}
provided $k_n \asymp\Delta_n^{-1/2}$. This assumption on the rate of
growth of $k_n$ guarantees that the distribution of $Z_i^n$ is
``sufficiently close'' to standard normal.
Now we can use~(\ref{eqproofevtt}) to get
%
%
\begin{eqnarray}
\mathbb{P}\bigl(c_{n}^{-1}(M_n-m_{n})
\leq x\bigr) &=& \bigl[\mathbb{P} \bigl(\bigl|Z_i^n\bigr|\leq
\tau_n(x) \bigr)\bigr]^{b_n}
\nonumber
\\[-8pt]
\\[-8pt]
\nonumber
&\sim&\bigl[\Phi\bigl(
\tau_n(x)\bigr)-\Phi\bigl(-\tau_n(x)\bigr)
\bigr]^{b_n}
\end{eqnarray}
as $n\rightarrow\infty$. From here, using the results for the maximum
domain of attraction of the Gumbel distribution (see e.g., Example~1.1.7 of~\cite{Haan}), we have
%
%
\begin{eqnarray}
\bigl[\Phi\bigl(\tau_n(x)\bigr)-\Phi\bigl(-\tau_n(x)
\bigr)\bigr]^{b_n} &=& \bigl[2\Phi\bigl(\tau _n(x)\bigr)-1
\bigr]^{b_n}
\nonumber
\\[-8pt]
\\[-8pt]
\nonumber
& \longrightarrow &\exp\bigl(-2\exp(-x)\bigr) \qquad \forall x,
\end{eqnarray}
and hence
%
%
\begin{equation}
c_{n}^{-1}(M_n-m_{n}) \stackrel{
\mathcal{L}} {\longrightarrow } \Lambda
\end{equation}
for $\Lambda$ being a random variable with c.d.f. $\exp(-2\exp(-x))$. From
here the result in (\ref{thmevt1}) (with $\eta_n$ replaced by $\eta
_n^*$) follows.

\textit{Step} 2. We now prove (\ref{thmevt1}) with condition (iii) relaxed.
Let $\eta_n^{\prime\ast}$ be defined as $\eta_n^*$ but with
$\widehat
{V}_t^{\ast}$ replaced by $\widehat{V}_t^{\prime}$. By step 1, (\ref
{thmevt1}) holds with $\eta_n$ replaced by $\eta_n^{\prime\ast}$. It
remains to show that
%
%
\begin{equation}
\log\bigl(\lfloor T/u_n\rfloor\bigr)^{1/2}
\Delta_n^{-1/4} \bigl(\eta_n - \eta
_n^{\prime\ast}\bigr) = o_p(1).\label{eta301}
\end{equation}
Note that $|\eta_n-\eta_n^{\prime\ast}| \leq\sup_{0\leq i \leq
\lfloor T/u_n \rfloor-1} |\widehat{V}_{iu_n} - \widehat{V}^{\prime
}_{iu_n}| = O_p(\Delta_n^{(2-r)\varpi-1/2})$,\vspace*{2pt} where the stochastic
order is shown in step 2 of the proof of Theorem~\ref{thm-u}. (\ref
{eta301}) then follows condition (ii). This completes the proof.
\end{pf*}

\subsection{\texorpdfstring{Proofs in Section \protect\ref{sec-rd}}{Proofs in Section 4.2}}

Under Assumption~\ref{assC}, by It\^{o}'s formula, we can represent $V$ as
\begin{eqnarray*}
V_{t} &=&V_{0}+\int_{0}^{t}
\tilde{b}_{V,s}\,ds+\int_{0}^{t}\tilde {
\sigma}%
_{V,s}\,dW_{s}+\int_{0}^{t}
\tilde{\sigma}_{V,s}^{\prime
}\,dW_{s}^{\prime}
\\
&&{}+\int_{0}^{t}\int_{\mathbb{R}}
\tilde{\delta}_{V} ( s,z ) ( \mu-\nu ) ( ds,dz ),
\end{eqnarray*}
where, by localization, we can assume without loss that the
coefficients $%
\tilde{b}_{V}$, $\tilde{\sigma}_{V}$, $\tilde{\sigma}_{V}^{\prime
}$ are
bounded, and $|\tilde{\delta}_{V} ( \omega,s,z ) |\leq
\tilde{
\Gamma} ( z ) $ for any $ ( \omega,s,z ) $, where
$\tilde{%
\Gamma} ( \cdot ) $ is bounded and deterministic, and
satisfies $%
\int_{\mathbb{R}}\tilde{\Gamma} ( z ) ^{\tilde
{r}}\lambda (
dz ) <\infty$.

We consider the following decomposition: for $q>0$,
\begin{eqnarray*}
%
%
V_{t}&=&V_{t}^{\prime} ( q ) +V_{t}^{\prime\prime}
( q ),\qquad\mbox{where }
\\
%
V_{t}^{c}&=&V_{0}+
\int_{0}^{t}\tilde{b}_{V,s}\,ds+\int
_{0}^{t}\tilde {\sigma} 
_{V,s}\,dW_{s}+\int_{0}^{t}
\tilde{\sigma}_{V,s}^{\prime
}\,dW_{s}^{\prime
},
\\
%
V_{t}^{\prime}
( q ) &=&V_{t}^{c}+\int_{0}^{t}
\int_{ \{
z\dvtx \tilde{%
\Gamma} ( z ) \leq q \} }\tilde{\delta}_{V} ( s,z ) ( \mu-\nu ) (
ds,dz )
\\
%
&&{}-\int
_{0}^{t}\int_{ \{
z\dvtx \tilde
{%
\Gamma} ( z ) >q \} }\tilde{
\delta}_{V} ( s,z ) \nu ( ds,dz ),
\\
%
V_{t}^{\prime\prime}
( q ) &=&\int_{0}^{t}\int_{ \{
z\dvtx \tilde{%
\Gamma} ( z ) >q \} }
\tilde{\delta}_{V} ( s,z ) \mu ( ds,dz ).%
\end{eqnarray*}
Denote $I ( n,i ) =[iu_{n}, ( i+1 ) u_{n})$. We also
set $%
\mathcal{I}_{n} ( q ) =\{0 \leq i\leq \lfloor
T/u_{n} \rfloor-1\dvtx \mu(I ( n,i ) \times\{z\dvtx \tilde
{\Gamma
}%
( z ) >q\})=0\}$ and $\mathcal{T}_{n} ( q )
=\bigcup_{i\in\mathcal{I}_{n} ( q ) }I ( n,i ) $.
Here, $%
\mathcal{I}_{n} ( q ) $ collects indices of intervals not
containing ``big'' jumps. We can
decompose $%
\widehat{F}_{n,T} ( x ) =\widehat{F}_{n,T} (
x;q ) +%
\widehat{R}_{n,T} ( x;q )$ where
\[
\widehat{F}_{n,T} ( x;q ) =u_{n}\sum
_{i\in\mathcal
{I}_{n} (
q ) }1_{ \{ \widehat{V}_{iu_{n}}\leq x \} },\qquad \widehat{R}_{n,T} ( x;q ) =
\int_{ [ 0,T ]
\setminus
\mathcal{T}_{n} ( q ) }1_{ \{ \widehat{V}_{s}\leq
x
\} }\,ds.%
\]
Analogously, we have $F_{T}(x)=F_{n,T}(x;q)+R_{n,T}(x;q)$, where
\[
F_{n,T} ( x;q ) =\int_{\mathcal{T}_{n} ( q )
}1_{
\{
V_{s}\leq x \} }\,ds,\qquad
R_{n,T} ( x;q ) =\int_{ [
0,T ] \setminus\mathcal{T}_{n} ( q ) }1_{ \{
V_{s}\leq
x \} }\,ds.
\]
Finally, we set
\begin{eqnarray*}
\hat{\eta}_{n} ( q ) &=&\sup_{i\in\mathcal{I}_{n} (
q )
}\llvert
\widehat{V}_{iu_{n}}-V_{iu_{n}}\rrvert,
\\
\eta_{n}^{\prime} ( q ) &=&\sup_{0\leq i\leq
\lfloor
T/u_{n} \rfloor}\sup
_{s\in\lbrack iu_{n}, ( i+1 )
u_{n})}\bigl\llvert V_{s}^{\prime} ( q )
-V_{iu_{n}}^{\prime
} ( q ) \bigr\rrvert,
\\
\eta_{n} ( q ) &=&\hat{\eta}_{n} ( q ) +\eta
_{n}^{\prime} ( q ).
\end{eqnarray*}

We now generalize Lemma~\ref{lem-fu} as follows.\vadjust{\goodbreak}

\begin{lemma}
\label{lem-rd}Suppose $\eta_{n} ( q_{n} ) =O_{p} (
w_{n} ) $ for some nonrandom sequences $q_{n}\rightarrow0$ and $%
w_{n}\rightarrow0$ and Assumption~\ref{assSA} with $r=2$. Then \textup{(a)} under
Assumption~\ref{assB}, $\widehat{F}_{n,T} (
x ) -F_{T} ( x ) =O_{p} ( w_{n} )
+O_{p} (
u_{n}q_{n}^{-\tilde{r}} ) $; \textup{(b)} under Assumption~\ref{assBprime}, the
assertion in
\textup{(a)} holds uniformly in $x\in\mathbb{R}$ and moreover \textup{(c)} $|\widehat
{Q}_{n,T}(\alpha)
-Q_{T}(\alpha)|\leq\xi_{n}\sup_{x\in\mathbb{R}%
}|\widehat{F}_{n,T}(x)-F_{T}(x)|$ for some tight sequence of variables
$\xi
_{n}$, provided $f_{T}(Q_{T}(\alpha))>0$ a.s.
\end{lemma}

\begin{pf}
(a) Observe that $\mathbb{E}[\int_{ [
0,T%
] \setminus\mathcal{T}_{n} ( q_{n} ) }\,ds]\leq
Ku_{n}q_{n}^{-%
\tilde{r}}$, yielding%
%
%
\begin{equation}
\qquad\sup_{x\in\mathbb{R}}R_{n,T} ( x;q_{n} )
=O_{p} \bigl( u_{n}q_{n}^{-%
\tilde{r}} \bigr), \qquad\sup_{x\in\mathbb{R}}\widehat{R}%
_{n,T} (
x;q_{n} ) =O_{p} \bigl( u_{n}q_{n}^{-\tilde
{r}}
\bigr). \label{lem-rd1-101}
\end{equation}
By definition, $\widehat{F}_{n,T} ( x;q_{n} ) =\int_{\mathcal
{T}%
_{n} ( q_{n} ) }1_{ \{ \widehat{V}_{s}\leq x \}
}\,ds$. Note
that over $\mathcal{T}_{n} ( q_{n} ) $, the process
$V_{t}^{\prime
\prime} ( q_{n} ) $ is identically zero. Hence, $\sup_{t\in
\mathcal{T}_{n} ( q_{n} ) }|\widehat{V}_{t}-V_{t}|\leq
\eta
_{n} ( q_{n} ) $. By a similar argument as in (\ref{lem-fu-101}),
we deduce%
%
%
\begin{eqnarray}\label{lem-rd1-201}
%
%
&&\bigl\llvert \widehat{F}_{n,T} ( x;q_{n} ) -F_{n,T}
( x;q_{n} ) \bigr\rrvert
\nonumber\\
%
&&\qquad\leq F_{n,T}
\bigl( x+\eta_{n} ( q_{n} );q_{n} \bigr)
-F_{n,T} \bigl( x-\eta_{n} ( q_{n} );q_{n} \bigr)
\\
%
&&\qquad\leq F_{T}
\bigl( x+\eta_{n} ( q_{n} ) \bigr) -F_{T} \bigl(
x-\eta_{n} ( q_{n} ) \bigr).\nonumber
\end{eqnarray}
By an argument similar to part (b) of Lemma~\ref{lem-fu}, we derive $%
F_{T} ( x+\eta_{n} ( q_{n} )  ) -F_{T} (
x-\eta
_{n} ( q_{n} )  ) =O_{p} ( w_{n} ) $. The assertion
of part (a) then follows (\ref{lem-rd1-101}) and~(\ref{lem-rd1-201}).

(b) By localization, we can suppose that $V$ is bounded and thus $%
f_{T} ( \cdot ) $ is compactly supported. Since $f_{T} (
\cdot
) $ is continuous, $\sup_{x\in\mathbb{R}}f_{T} (
x ) =O_{p} ( 1 ) $. By~(\ref{lem-rd1-201}),%
%
%
\begin{equation}
\sup_{x\in\mathbb{R}}\bigl|\widehat{F}_{n,T} ( x;q_{n} )
-F_{n,T} ( x;q_{n} ) \bigr|\leq2\eta_{n} (
q_{n} ) \sup_{z\in
\mathbb{R}
}f_{T} ( z ).
\label{lem-rd-203}
\end{equation}
The assertion then readily follows (\ref{lem-rd1-101}) and (\ref
{lem-rd-203}%
).

(c) Observe that $\widehat{F}_{n,T}(\widehat{Q}_{n,T}(\alpha))\geq
\alpha
=F_{T}(Q_{T}(\alpha))$, where the inequality follows the definition of
quantiles and the equality is due to the continuity of $F_{T}(\cdot)$.
Hence,%
%
%
\begin{eqnarray}\label{lem-rd-304}
%
%
F_{T}\bigl(Q_{T}(\alpha)\bigr)-F_{T}\bigl(
\widehat{Q}_{n,T}(\alpha)\bigr) & \leq& 
%
\widehat{F}_{n,T}\bigl(
\widehat{Q}_{n,T}(\alpha)\bigr)-F_{T}\bigl(\widehat
{Q}_{n,T}(\alpha )\bigr)
\nonumber
\\[-8pt]
\\[-8pt]
\nonumber
%
 & \leq&
%
\sup_{x\in\mathbb{R}}
\bigl\llvert \widehat {F}_{n,T}(x)-F_{T}(x)\bigr\rrvert.%
\end{eqnarray}
For any $\varepsilon>0$, $\widehat{F}_{n,T}(\widehat{Q}_{n,T}(\alpha
)-\varepsilon)<\alpha=F_{T}(Q_{T}(\alpha))$, yielding
\begin{eqnarray*}
F_{T}\bigl(\widehat{Q}_{n,T}(\alpha)-\varepsilon
\bigr)-F_{T}\bigl(Q_{T}(\alpha)\bigr) &<&F_{T}
\bigl(\widehat{Q}_{n,T}(\alpha)-\varepsilon\bigr)-\widehat
{F}_{n,T}\bigl(\widehat{%
Q}_{n,T}(\alpha)-
\varepsilon\bigr)
\\
&\leq&\sup_{x\in\mathbb{R}}\bigl\vert \widehat {F}_{n,T}(x)-F_{T}(x)
\bigr\vert.
\end{eqnarray*}
Since $F_{T} ( \cdot ) $ is continuous, by sending
$\varepsilon
\downarrow0$ we deduce%
%
%
\begin{equation}
F_{T}\bigl(\widehat{Q}_{n,T}(\alpha)\bigr)-F_{T}
\bigl(Q_{T}(\alpha)\bigr)\leq\sup_{x\in
\mathbb{R}}\bigl\llvert
\widehat{F}_{n,T}(x)-F_{T}(x)\bigr\rrvert.
\label{lem-rd-305}
\end{equation}
Let $\mathcal{K}_{n,T} ( \alpha ) $ be the closed interval with
endpoints $Q_{T}(\alpha)$ and $\widehat{Q}_{n,T}(\alpha)$. Set $\xi
_{n}\equiv\sup_{x\in\mathcal{K}_{n,T} ( \alpha )
}f_{T}^{-1} ( x ) $. By a mean-value expansion,%
%
%
\begin{equation}
\qquad\bigl\llvert F_{T}\bigl(\widehat{Q}_{n,T}(\alpha)
\bigr)-F_{T}\bigl(Q_{T}(\alpha )\bigr)\bigr\rrvert \geq\inf
_{x\in\mathcal{K}_{n,T} ( \alpha )
}f_{T} ( x ) \bigl\llvert \widehat{Q}_{n,T}(
\alpha )-Q_{T}(\alpha )\bigr\rrvert. \label{lem-rd-306}
\end{equation}
Since $f_{T}(Q_{T}(\alpha))>0$ a.s., $\widehat{Q}_{n,T}(\alpha
)\stackrel{%
\mathbb{P}}{\longrightarrow}Q_{T}(\alpha)$ by Corollary~\ref
{cor-Q}. Since
$f_{T} ( \cdot ) $ is continuous, $\inf_{x\in\mathcal{K}%
_{n,T} ( \alpha ) }f_{T} ( x ) \stackrel
{\mathbb{P}}{%
\longrightarrow}f_{T}(Q_{T}(\alpha))>0$; hence $\xi_{n}$ is tight. The
assertion then follows (\ref{lem-rd-304}), (\ref{lem-rd-305}) and
(\ref%
{lem-rd-306}).
\end{pf}

\begin{pf*}{Proof of Theorem~\ref{thm-rd}} \textit{Step} 1. We first
consider $%
\eta_{n}^{\prime} ( q_{n} ) $. For each $i$,
\begin{eqnarray*}
%
%
&&\sup_{s\in\lbrack iu_{n}, ( i+1 ) u_{n})}\bigl\llvert V_{s}^{\prime
} (
q_{n} ) -V_{iu_{n}}^{\prime} ( q_{n} ) \bigr
\rrvert \leq\zeta_{n,i}+\zeta_{n,i}^{\prime}+
\zeta_{n,i}^{\prime\prime
},\qquad\mbox{where}
\\
%
&&\qquad\zeta_{n,i}=\sup
_{s\in\lbrack iu_{n}, ( i+1 )
u_{n})}\biggl\llvert \int_{iu_{n}}^{s}
\int_{ \{ z\dvtx \tilde{\Gamma} ( z )
>q_{n} \} }\tilde{\delta}_{V} ( s,z ) \nu ( ds,dz )
\biggr\rrvert,
\\
%
&&\qquad\zeta_{n,i}^{\prime}=
\sup_{s\in\lbrack iu_{n}, ( i+1 )
u_{n})}\biggl\llvert \int_{iu_{n}}^{s}
\int_{ \{ z\dvtx \tilde{\Gamma
} (
z ) \leq q_{n} \} }\tilde{\delta}_{V} ( s,z ) ( \mu -\nu ) (
ds,dz ) \biggr\rrvert,
\\
%
&&\qquad\zeta_{n,i}^{\prime\prime}=
\sup_{s\in\lbrack iu_{n}, (
i+1 )
u_{n})}\bigl\llvert V_{s}^{c}-V_{iu_{n}}^{c}
\bigr\rrvert.
\end{eqnarray*}
For $\zeta_{n,i}$, observe that%
%
%
\begin{eqnarray}\label{thm-rd-101}\quad
%
%
\Bigl\llvert \sup_{0\leq i\leq \lfloor T/u_{n} \rfloor
}\zeta _{n,i}\Bigr\rrvert &
\leq& 
%
\sup_{0\leq i\leq \lfloor T/u_{n} \rfloor}
\int_{iu_{n}}^{ (
i+1 ) u_{n}}\int_{ \{ z\dvtx \tilde{\Gamma} ( z )
>q_{n} \} }\tilde{
\Gamma} ( z ) \nu ( ds,dz )
\nonumber
\\[-8pt]
\\[-8pt]
\nonumber
\quad & \leq& 
%
Ku_{n}q_{n}^{ ( 1-\tilde{r} ) \wedge0}.
\end{eqnarray}

Now turn to $\zeta_{n,i}^{\prime}$. Let $p\geq2$. For each $i\geq
0$, by
Lemma~2.1.5 in~\cite{jacodprotter2012},
\begin{eqnarray*}
\mathbb{E}\bigl\llvert \zeta_{n,i}^{\prime}\bigr\rrvert
^{p} &\leq &K_{p}u_{n}\int_{ \{ z\dvtx \tilde{\Gamma} ( z ) \leq
q_{n} \} }%
\tilde{\Gamma} ( z ) ^{p}\lambda ( dz )
\\
&&{}+K_{p}u_{n}^{p/2} \biggl( \frac{1}{u_{n}}\int
_{I ( n,i )
}\,ds\int_{ \{ z\dvtx \tilde{\Gamma} ( z ) \leq q_{n}
\}
}\tilde{%
\Gamma} ( z ) ^{2}\lambda ( dz ) \biggr) ^{p/2}
\\
&\leq&K_{p}u_{n}q_{n}^{p-\tilde{r}}+K_{p}u_{n}^{p/2}q_{n}^{ (
2-\tilde{r%
} ) p/2}.
\end{eqnarray*}
Hence, $\Vert\zeta_{n,i}^{\prime}\Vert_{p}\leq
K_{p}u_{n}^{1/p}q_{n}^{1-\tilde{r}/p}+K_{p}u_{n}^{1/2}q_{n}^{1-\tilde{r}/2}$.
By the maximal inequality (Lemma~2.2.2 in~\cite{vandervaartwellner}),%
%
%
\begin{equation}
\Bigl\llVert \sup_{0\leq i\leq \lfloor T/u_{n} \rfloor
}\zeta _{n,i}^{\prime}
\Bigr\rrVert _{p}\leq K_{p}q_{n}^{1-\tilde{r}%
/p}+K_{p}u_{n}^{1/2-1/p}q_{n}^{1-\tilde{r}/2}.
\label{thm-rd-102}
\end{equation}
Since $V^{c}$ is continuous, by a standard estimate for $\zeta
_{n,i}^{\prime\prime}$ and the maximal inequality,%
%
%
\begin{equation}
\Bigl\llVert \sup_{0\leq i\leq \lfloor T/u_{n} \rfloor
}\zeta _{n,i}^{\prime\prime}
\Bigr\rrVert _{p}\leq K_{p}u_{n}^{1/2-1/p}.
\label{thm-rd-103}
\end{equation}
Combining (\ref{thm-rd-101}), (\ref{thm-rd-102}) and (\ref
{thm-rd-103}), we
derive for $p\geq2$, $\Vert\eta_{n}^{\prime}( q_{n}) \Vert_{p}
\leq
K_{p}a_{n,p}^{\prime}$, where
\[
a_{n,p}^{\prime} \equiv u_{n}q_{n}^{ ( 1-\tilde{r} )
\wedge0}
\vee q_{n}^{1-\tilde{r}/p}\vee u_{n}^{1/2-1/p}.
\]

\textit{Step} 2. Observe that%
%
%
\begin{equation}
\hat{\eta}_{n} ( q_{n} ) \leq\sup_{i\in\mathcal
{I}_{n} (
q_{n} ) }
\bigl\llvert \widehat{V}_{iu_{n}}-\widehat {V}_{iu_{n}}^{\prime
}
\bigr\rrvert +\sup_{i\in\mathcal{I}_{n} ( q_{n} ) }\bigl\llvert \widehat{V}_{iu_{n}}^{\prime}-V_{iu_{n}}
\bigr\rrvert. \label{thm-rd-201}
\end{equation}
By Lemma~\ref{lem-vstar}, for $i\in\mathcal{I}_{n} (
q_{n} )
$, for some $\xi_{n,i}$ with $\Vert\xi_{n,i} \Vert_p \leq K_p k_n^{-1/2}$,
\[
\bigl\llvert \widehat{V}_{iu_{n}}^{\prime}-V_{iu_{n}}\bigr
\rrvert \leq \xi_{n,i}+\sup_{s\in\lbrack iu_{n}, ( i+1 )
u_{n})}\bigl\llvert
V_{s}^{\prime} ( q_{n} ) -V_{iu_{n}}^{\prime
}
( q_{n} ) \bigr\rrvert.
\]
Therefore, by a similar argument as in step 1,%
%
%
\begin{equation}
\Bigl\llVert \sup_{i\in\mathcal{I}_{n} ( q_{n} ) }\bigl\llvert \widehat{V}_{iu_{n}}^{\prime}-V_{iu_{n}}
\bigr\rrvert \Bigr\rrVert _{p}\leq K_{p}u_{n}^{-1/p}k_{n}^{-1/2}+K_{p}a_{n,p}^{\prime}.
\label{thm-rd-202}
\end{equation}
Similarly as in step 2 of the proof of Theorem~\ref{thm-u} [recall
that $\theta
=0$ if $r\leq1$ and $\theta\in ( 0,\infty ) $ can be
arbitrarily fixed when $r>1$]%
%
%
\begin{equation}
\sup_{i\in\mathcal{I}_{n} ( q_{n} ) }\bigl\llvert \widehat {V}%
_{iu_{n}}-\widehat{V}_{iu_{n}}^{\prime}\bigr\rrvert
=O_{p} \bigl( \Delta _{n}^{{(\gamma-r\varpi-\theta)}/{(1\vee r)}- ( 1-2\varpi
)
} \bigr).
\label{thm-rd-203}
\end{equation}
Combining (\ref{thm-rd-201})--(\ref{thm-rd-203}), we derive $\hat
{\eta}%
_{n} ( q_{n} ) =O_{p} ( w_{n,p}  ) $ for $%
p\geq2$, where
\[
w_{n,p} \equiv\Delta_{n}^{{(\gamma-r\varpi-\theta)}/{(1\vee r)}%
- ( 1-2\varpi ) }\vee
u_{n}^{-1/p}k_{n}^{-1/2}\vee
a_{n,p}^{\prime
}.
\]
By step 1, we further derive $\eta_{n} ( q_{n} )
=O_{p} (
w_{n,p} ) $.

By Lemma~\ref{lem-rd}(a), we have
\begin{eqnarray*}
%
%
&&\widehat{F}_{n,T} ( x ) -F_{T} ( x )
\\
%
&&\qquad=O_{p} (
w_{n,p} ) +O_{p} \bigl( u_{n}q_{n}^{-%
\tilde{r}}
\bigr)
\\
%
&&\qquad=O_{p} \bigl(
\Delta_{n}^{{(\gamma-r\varpi-\theta)}/{(1\vee r)}
- ( 1-2\varpi ) }\vee u_{n}^{-1/p}k_{n}^{-1/2}
\vee q_{n}^{1-\tilde{%
r}/p}\vee u_{n}^{1/2-1/p}\\
&&\hspace*{280pt}{}\vee
u_{n}q_{n}^{-\tilde{r}} \bigr).
\end{eqnarray*}
Taking $q_{n}=u_{n}^{1/ ( 1+\tilde{r}-\tilde{r}/p ) }$ and
recalling $k_{n}\asymp\Delta_{n}^{-\gamma}$, we have
\begin{eqnarray*}
%
%
&&\widehat{F}_{n,T} ( x ) -F_{T} ( x )
\\
%
&&\qquad=O_{p} \bigl(
\Delta_{n}^{{(\gamma-r\varpi-\theta)}/{(1\vee r)}
- ( 1-2\varpi ) }\vee\Delta_{n}^{{\gamma}/{2}-
{(1-\gamma)}/{p}}\vee
\Delta_{n}^{ ( 1-\gamma ) (
{1}/{2}-{1}/{p})}\\
&&\hspace*{193pt}{}\vee\Delta _{n}^{{((1-\gamma)(1-\tilde{r}/p))}/{(1+\tilde{r}-\tilde
{r}/p)}}
\bigr)
\end{eqnarray*}
and by taking $p$ sufficiently large,
\begin{eqnarray*}
%
%
&&\widehat{F}_{n,T} ( x ) -F_{T} ( x )
\\
%
&&\qquad=O_{p} \bigl(
\Delta_{n}^{{(\gamma-r\varpi-\theta)}/{(1\vee r)}
- ( 1-2\varpi ) }\vee\Delta_{n}^{\gamma/2-\iota}\vee
\Delta _{n}^{ ( 1-\gamma ) /2-\iota}\vee\Delta_{n}^{ (
1-\gamma
) / ( 1+\tilde{r} ) -\iota}
\bigr).
\end{eqnarray*}

The discontinuous case in part (a) then readily follows the definitions
of~$\theta$, $a_{n}$ [see (\ref{an})] and $d_{n}$. The continuous case in part
(a), as well as part (b), can be proved in a similar (but simpler)
way.
\end{pf*}

\begin{pf*}{Proof of Theorem~\ref{thm-rdu}} We first show that
$%
\sup_{x\in\mathbb{R}}|\widehat{F}_{n,T}(x)-F_{T}(x)|=O_{p}(d_{n})$. The
proof is similar as that of Theorem~\ref{thm-rd}, except in step 2 of the
proof, we use Lemma~\ref{lem-rd}(b) instead of Lemma~\ref{lem-rd}(a). The
result for $\widehat{F}_{n,T}^{\ast}$ can be proved similarly. The
assertion concerning $\widehat{Q}_{n,T}(\alpha)$ then follows from
Lemma~\ref{lem-rd}%
(c). The assertion on $\widehat{Q}^{*}_{n,T}(\alpha)$ can be proved in
a similar (but simpler) way; the details are omitted for brevity.
\end{pf*}

\subsection{\texorpdfstring{Proofs in Section \protect\ref{sec-ker}}{Proofs in Section 5}}
\mbox{}
\begin{pf*}{Proof of Theorem~\ref{thm-kd}} (a) By localization and
condition (iv), we can assume that $V_{t}$ takes value in some compact $
\mathcal{K}\subset ( 0,\infty ) $, and thus $f_{T} (
\cdot
) $ is supported on~$\mathcal{K}$. We set $f_{n,T}(x)=%
\int_{0}^{T}h_{n}^{-1}\kappa( h_{n}^{-1}(V_{s}-x)) \,ds$. For each
$x\in\mathbb{R}$,
\begin{eqnarray*}
\mathbb{E}\bigl\vert \widehat{f}_{n,T} ( x ) -f_{n,T} ( x
) \bigr\vert &\leq&\mathbb{E} \biggl[ h_{n}^{-1}\int
_{0}^{T}%
\biggl\vert\kappa 
%
\biggl(\frac{\widehat{V}_{s}-x}{h_{n}}%
%
\biggr)-\kappa 
%
\biggl(\frac{V_{s}-x}{h_{n}}%
%
\biggr)\biggr\vert \,ds \biggr]
\\
&\leq&Kh_{n}^{-2}\mathbb{E} \biggl[ \int
_{0}^{T}\vert \widehat {V}%
_{s}-V_{s}\rrvert ds \biggr].
\end{eqnarray*}
By Lemmas~\ref{lem-vstar} and~\ref{lem-eoj}, $\mathbb{E}|\widehat{V}
_{s}-V_{s}|\leq K\bar{a}_{n}$. Hence,
%
%
\begin{equation}
\mathbb{E}\bigl\llvert \widehat{f}_{n,T} ( x ) -f_{n,T} ( x
) \bigr\rrvert \leq Kh_{n}^{-2}\bar{a}_{n},\qquad
\mathbb {E} \bigl[ \Vert\widehat{f}_{n,T}-f_{n,T}
\Vert_{w} \bigr] \leq Kh_{n}^{-2}\bar{a}_{n},
\label{ker-1}\hspace*{-35pt}
\end{equation}
where $K$ does not depend on $x$. Now observe that $f_{n,T}(x)=\int_{\mathbb{%
R}}h_{n}^{-1}\kappa( h_{n}^{-1}(y-x)) f_{T}(y)\,dy$. By a change
of variable, $f_{n,T}(x)=\int_{\mathbb{R}}\kappa ( z )
f_{T}(x+h_{n}z)\,dz$. Hence,
\begin{eqnarray*}
\mathbb{E}\bigl\llvert f_{n,T}(x)-f_{T} ( x ) \bigr\rrvert
&\leq &\int_{\mathbb{R}}\kappa ( z ) \mathbb{E}\bigl\llvert
f_{T} ( x+h_{n}z ) -f_{T} ( x ) \bigr\rrvert \,dz
\\
&\leq&Kh_{n}^{\beta}\int_{\mathbb{R}}\kappa ( z )
\llvert z\rrvert ^{\beta}\,dz\leq Kh_{n}^{\beta},
\end{eqnarray*}
which further implies $\mathbb{E}[\Vert f_{n,T}-f_{T}\Vert_{w}]\leq
Kh_{n}^{\beta}$. Combining these estimates with~(\ref{ker-1})
completes the
proof of part (a). Part (b) can be proved in a similar way.
\end{pf*}

\section*{Acknowledgments}
We would like to thank Tim Bollerslev, Nathalie Eisenbaum, Jean Jacod,
Andrew Patton and Philip Protter for helpful discussions as well as an
Associate Editor and two referees for very helpful suggestions. We are
particularly grateful to Markus Reiss for suggesting the direct
estimation approach adopted in the paper and the link with the uniform
error in estimating the volatility path given in Lemma~\ref{lem-fu} of the paper.

%

%



\printaddresses

\end{document}